\documentclass[11pt,reqno]{amsart}
\usepackage[colorlinks=true,urlcolor=blue, citecolor=red,linkcolor=blue,linktocpage,pdfpagelabels, bookmarksnumbered,bookmarksopen]{hyperref}
\usepackage{xcolor}
\usepackage{amsthm} 
\usepackage{latexsym,amsmath,amssymb}
\usepackage{accents}
\usepackage[colorinlistoftodos,prependcaption,textsize=tiny]{todonotes}
\usepackage{a4wide}
\usepackage{soul}
\usepackage{mathtools} 
\usepackage{xparse} 
\usepackage{comment}
\usepackage{enumitem}
\usepackage{cancel}
\usepackage[sort, numbers]{natbib}

\title[Stationary systems with signed Coulombian interactions]{Mean-field limit of 2D stationary particle systems with signed Coulombian interactions}

\date{\today}

\author{Jan Peszek}
\address[J. Peszek ]{Institute of Mathematics, University of Warsaw, Banacha 2, 02-097
Warszawa, Poland}
\email{j.peszek@mimuw.edu.pl}

\author{Rémy Rodiac}
\address[R. Rodiac]{Institute of Mathematics, University of Warsaw, Banacha 2, 02-097
Warszawa, Poland}
\email{rrodiac@mimuw.edu.pl}

\definecolor{indigo}{rgb}{0.29, 0.0, 0.51}

\newcommand{\mres}{\mathbin{\vrule height 1.6ex depth 0pt width
0.13ex\vrule height 0.13ex depth 0pt width 1.3ex}}

plus 6pt minus 12pt
\def\Om{\Omega}

\def\C{{\mathcal C}}
\def\D{{\mathcal D}}

\def \O{{\Omega}}
\def\N{{\mathbb N}}

\def\H{{\mathcal H}}

\def\M{\mathcal{M}}

\def\e {{\varepsilon}}

\newcommand{\dd}{\,\mathrm{d}}

\newtheorem{theorem}{Theorem}
\newtheorem{lemma}[theorem]{Lemma}
\newtheorem{corollary}[theorem]{Corollary}
\newtheorem{proposition}[theorem]{Proposition}

\newtheorem{remark}[theorem]{Remark}
\newtheorem{definition}[theorem]{Definition}

\newcommand{\loc}{\mathrm{loc}}

\def\supp{{\rm supp\,}}
\def\loc{{\rm loc}}
\def \p {\partial}

\def\H{\mathcal{H}}
\def\capp{{\rm cap}}

\newcommand{\R}{\mathbb{R}}

\newcommand{\barint}{
\rule[.036in]{.12in}{.009in}\kern-.16in \displaystyle\int }

\newcommand{\barcal}{\mbox{$ \rule[.036in]{.11in}{.007in}\kern-.128in\int $}}

\def\mvint_#1{\mathchoice
          {\mathop{\vrule width 6pt height 3 pt depth -2.5pt
                  \kern -8pt \intop}\nolimits_{\kern -3pt #1}}%
          {\mathop{\vrule width 5pt height 3 pt depth -2.6pt
                  \kern -6pt \intop}\nolimits_{#1}}%
          {\mathop{\vrule width 5pt height 3 pt depth -2.6pt
                  \kern -6pt \intop}\nolimits_{#1}}%
          {\mathop{\vrule width 5pt height 3 pt depth -2.6pt
                  \kern -6pt \intop}\nolimits_{#1}}}

\numberwithin{theorem}{section} \numberwithin{equation}{section}

\newcommand{\lh}{[}
\newcommand{\rh}{]}
\newcommand{\vx}{{\bf x}}

\DeclareMathOperator \tr{\text{tr}}
\DeclareMathOperator \dive{{\rm div}}

\renewcommand\phi{\varphi}
\renewcommand\epsilon{\varepsilon}

\let\latexchi\chi
\makeatletter
\renewcommand\chi{\@ifnextchar_\sub@chi\latexchi}
\newcommand{\sub@chi}[2]{
  \@ifnextchar^{\subsup@chi{#2}}{\latexchi^{}_{#2}}%
}
\newcommand{\subsup@chi}[3]{
  \latexchi_{#1}^{#3}%
}

\begin{document}

\begin{abstract}
We study the mean-field limits of critical points of interaction energies with Coulombian singularity. An important feature of our setting is that we allow interaction between particles of opposite signs. Particles of opposite signs attract each other whereas particles of the same signs repel each other. 
In 2D, we prove that the associated empirical measures converge to a limiting measure \(\mu\) that satisfies a two-fold criticality condition: in velocity form or in vorticity form. Our setting includes the stationary attraction-repulsion problem with Coulombian singularity and the stationary system of point-vortices in fluid mechanics. In this last context, in the case where the limiting measure is in \(H^{-1}_{\text{loc}}(\R^2)\), we recover the classical criticality condition stating that \(\nabla^\perp g \ast \mu\), with \(g(x)=-\log |x|\), is a stationary solution of the incompressible Euler equation. This result, is, to the best of our knowledge, new in the case of particles with different signs (for particles of the positive sign it was obtained by Schochet in 1996).
In order to derive the limiting criticality condition in the velocity form, we follow an approach devised by Sandier-Serfaty in the context of Ginzburg-Landau vortices. This consists of passing to the limit in the stress-energy tensor associated with the velocity field. On the other hand, the criticality condition in the vorticity form is obtained by arguments closer to the ones of Schochet.
\end{abstract}

\maketitle

\section{Introduction}

We study the mean-field limits of solutions  \(\vx_N:= (x_{1},\dots,x_{N})\in (\R^d)^N\) to the equilibrium system
\begin{multline}\label{eq:critical_points}
\frac{1}{M_N}\sum_{j\neq i} d_{j} \left( \nabla g(x_{i}-x_{j})+\nabla F(x_{i}-x_{j})\right) +  \nabla V(x_{i})=0, 
\quad x_{i}\neq x_{j},\ \forall i\neq j\in\{1,\dots,N\}.
\end{multline}
Here \(d_{j}\) are real mass (or charge) coefficients, \(M_N:=\sum_{i=1}^N |d_{i}|\), and the involved functions $g + F$ and $V$ represent interaction and potential energy, respectively. 
The key features of our setting are its Coulombian singularity (expressed in $g$) and the admissibility of negative topological charges $d_i$. That is to say, throughout the paper, we assume that
\begin{equation*}
g(x)= \begin{cases}
-\log |x| & \text{for } d=1,2, \\
\frac{1}{|x|^{d-2}} & \text{for } d\geq 3,
\end{cases}\quad
-\Delta g= c_d \delta_0, \quad c_d=\begin{cases}
    2\pi & \text{ if } d=2 \\
    (d-2)\mathcal{H}^{d-1}(\mathbb{S}^{d-1}) & \text{ if } d\geq 3
\end{cases}
\end{equation*}
is, up to a constant, the fundamental solution of the Laplace equation.
Taking into account the regular part  $F$ of the interaction energy and varying signs of $d_i$, we note that \eqref{eq:critical_points} describes, in particular, the stationary case of the  attraction-repulsion system of interacting agents, cf.\ e.g.\  \cite{Carrillo_Choi_Perez_2017, C-W-Z_2020}. This model is also used to describe stationary systems of point vortices  or relative equilibria\footnote{ The point-vortex system reads \(x_i'(t)=\sum_{j\neq i} d_j \nabla^\perp g(x_i(t)-x_j(t))\). Stationary solutions are thus solutions of \eqref{eq:critical_points} with \(F=V=0\). Relative equilibria can be of the form \(x_i(t)=x_i(0)+\alpha t\) for some \(\alpha=(\alpha_1,\alpha_2)\in \R^2\). In that case they should satisfy \(\alpha=\frac{1}{M_N}\sum_{j\neq i} d_j\nabla^\perp g(x_i(0)-x_j(0))=0\) for all \(i=1,\dots,N\) which is equivalent to \eqref{eq:critical_points} with \(V(x)=-\alpha_2 x_1+\alpha_1 x_2\). Another type of relative equilibrium is of the form \(x_i(t)= x_i(0) e^{i\omega t}\), for some \(\omega \in \R\). In that case the equation is \(i\omega x_i(0)=\frac{1}{M_N}\sum_{j\neq i} d_j \nabla^\perp g (x_i(0)-x_j(0))\) and is equivalent to \eqref{eq:critical_points} with \(V(x)=\omega |x|^2\). } cf.\ \cite{Aref_2007}.

Solutions $\vx_N$ of \eqref{eq:critical_points} are, by definition, critical points of the Hamiltonian
\begin{equation*}
\H_N(\vx)=\frac{1}{2M_N}\sum_{i=1}^N \sum_{j\neq i} d_id_j\left( g(x_i-x_j)+F(x_i-x_j)\right) +\sum_{i=1}^N d_iV(x_i)
\end{equation*}
and their mean-field limits, i.e. weak limits $\mu$ of empirical measures $\mu_N:=\frac{1}{M_N}\sum_{i=1}^N d_i\delta_{x_i}$, are expected to satisfy certain criticality conditions for the continuum energy
\begin{equation}\label{eq:continuum_energy_2}
    \mathcal{H}(\mu):= \frac12 \iint_{\R^d\times \R^d} \left(g(x-y)+F(x-y)\right) \dd \mu(x)\dd \mu(y) +\int_{\R^d}V(x) \dd \mu(x).
\end{equation}
 In fact, our goal is to justify the link between the discrete energy \(\H_N\) and the continuum energy \(\H\) by showing that, in dimension 2 and under some assumptions on \(F\) and \(V\),


\smallskip
\begin{center}
    {\it critical points of  $\H_N$ converge to critical points of  $\H$,}\\

    \medskip
    {\it stable critical points of $\H_N$ converge to stable critical points of $\H$.}
\end{center}
\medskip

An essential part of our investigation is to find and discuss various suitable notions of criticality for the continuum energy $\H$. Thus, the rigorous formulation of the main results is delicate and requires some preliminary information --  it can be found in Section \ref{sec:results} below. Meanwhile, let us explain the origins of our approach.  We follow two different main ideas: one due to Sandier-Serfaty \cite{Sandier_Serfaty_2007} and the other due to Schochet \cite{Schochet_1995}.

For notational simplicity hereinafter we write \( (x_1,\dots,x_N)\) and \(d_i\), even though technically  it should be \( (x_{1,N},\dots,x_{N,N})\) and \(d_{i,N}\), due to the dependence on $N$. We look at  critical points of $\H_N$, i.e. solutions to \eqref{eq:critical_points}, from three distinct points of view
\begin{equation}\label{eq:3crit}
\begin{split}
    \vx_N &= (x_1,...,x_N)\in(\R^d)^N,\quad \mu_N:=\frac{1}{M_N}\sum_{i=1}^N d_i\delta_{x_i} \quad \text{with }
    M_N = \sum_{i=1}^N |d_{i}|,\\
    h_N &=g \ast \mu_N=\frac{1}{M_N}\sum_{i=1}^N d_ig(\cdot-x_i),
\end{split}
\end{equation}
where $\vx_N$ is a classical critical point of $\H_N$, the measure $\mu_N$ is its associated {\it empirical measure}, while $h_N$ is its {\it electric potential}.

For the first approach towards the notion of criticality for the continuous energy $\H$, we employ the electric potential $h_N$.
We remark that Equation \eqref{eq:critical_points} means that for the system of interacting points to be at equilibrium, the force felt by the point \(x_i\) and generated by all the other points \(x_j\) with \(j \neq i\) should be zero. This formally reads
\begin{equation}\label{eq:formal_equilibrium_1}
\mu_N \nabla (h_N+F\ast \mu_N+V)=0.
\end{equation}
Formally again, passing to the limit in \eqref{eq:formal_equilibrium_1} leads to the conclusion that \(\mu\) satisfies the critical condition \(\mu  (\nabla g \ast \mu+\nabla F\ast \mu+V)=0\) in \(\R^d\). In the first part of this article, we use the stress-energy tensor associated with the Dirichlet energy and the notion of divergence-free in finite part vector fields, due to Sandier-Serfaty \cite{Sandier_Serfaty_2007}, to give a rigorous formulation of \eqref{eq:formal_equilibrium_1} and to pass to the limit in this relation. 

In the second approach, inspired by Schochet \cite{Schochet_1995}, we adopt another point of view. We show that $\vx_N$ satisfying \eqref{eq:critical_points} verify
\begin{multline}\label{eq:vorticity_form}
    \frac12 \iint_{(\R^2 \times \R^2) \setminus \{x=y\}}\left( \nabla g(x-y)+\nabla F(x-y)\right)\cdot   \left(\varphi(x)-\varphi(y)\right)\dd \mu_N(x) \dd \mu_N(y) \\
    + \int_{\R^2}\nabla V(x)\cdot \varphi(x) \dd \mu_N(x)=0,
\end{multline}
for all \(\varphi \in \C^\infty_c(\R^2;\R^2)\), and we pass to the limit as $N\to +\infty$ in this relation.

We conclude the first part of the introduction by stating our assumptions on \(F\) and \(V\) appearing in \(\H_N\) and our assumptions on the sequences of critical points we consider. We always assume that 
$$ F\in \C^1(\R^d)\quad \mbox{ and }\quad F(x)=F(-x) \mbox{ for all } x\in\R^d, $$
and that
\begin{align}
  V\in \C^1(\R^d) \text{  is bounded from below and }   \lim_{|x| \to  +\infty} \left( \frac{V(x)}{2}+g(x) \right)= +\infty. \label{A2}
\end{align}

\noindent
Then, to ensure the tightness of the sequence of empirical measures $(\mu_N)_N$, we assume that

\begin{equation}\label{eq:Boundedness_Hamiltonian}
    \sup_N \frac{1}{2M_N^2}\sum_{i=1}^N \sum_{j\neq i} |d_i||d_j|\left( g(x_i-x_j)+|F(x_i-x_j)|\right) +\frac{1}{M_N} \sum_{i=1}^N |d_i| V(x_i)< +\infty.
\end{equation}

\noindent
When \((\mu_N)_N\) is a sequence of finite signed Radon measures, satisfying the boundedness condition \eqref{eq:Boundedness_Hamiltonian}, it can be shown, cf.\ Lemma \ref{lem:tightness} below, that
there exists a subsequence (not relabelled) and a limit Radon measure \(\mu\) such that \(\mu_N \to \mu\) narrowly (tested by bounded continuous functions). We discuss other admissible assumptions, leading to stronger types of uniform integrability, later in Remark \ref{rem:otherass}.

\subsection{Previous results}

 When considering evolution problems, one can ask if gradient flows or Hamiltonian flows for discrete energies of the form $\H_N$ converge to their corresponding flows for the continuum energies $\H$. When all the charges of the particles are positive this question was studied e.g.\ in \cite{Goodman_Hou_Lowengrub_1990,Schochet_1996,Serfaty_2020,Nguyen_Rosenzweig_Serfaty_2022,Rosenzweig_2022,Rosenzweig_Serfaty_2023,Bresch_Jabin_Wang_2023}.  The case of signed charges was studied by Hauray \cite{Hauray_2009}; he considered the mean-field limits of evolving interacting points of Euler-type equations with the  kernel \(g\) having a strictly stronger singularity than the Coulombian one. In 1D, results for a variety of types of interactions were obtained in \cite{Carrillo_Ferreira_Precioso_2012} and \cite{Berman_Onnheim_2019}.

When all the \(d_i\)'s are positive, under the assumptions \(F=0\) and \eqref{A2},
it is known that the empirical measures associated with minimisers of $\H_N$ converge to minimisers of 
\begin{equation}\label{eq:continuum_energy}
   \frac12 \iint_{\R^d\times \R^d} g(x-y) \dd \mu(x)\dd \mu(y) +\int_{\R^d}V(x) \dd \mu(x),
\end{equation}
defined on the set of probability measures \(\mathcal{P}(\R^d)\) cf.\ \cite[Theorem 2.2]{Serfaty_2015}.   We also refer to  \cite[Theorem 1.2]{Canizo_Patacchini_2018} for convergence of minimisers of discrete sub-Coulombian energies with \(F\neq 0\), \(V=0\) and \(d_i\geq 0\), to minimisers of a continuum energy of the form \eqref{eq:continuum_energy_2}.

Assuming that \eqref{A2} holds true, that \(F=0\) and \(\{x\in \R^d: V(x)< +\infty\}\) has positive capacity, a minimiser of \eqref{eq:continuum_energy} exists, is unique (since in that case the energy \eqref{eq:continuum_energy} is strictly convex) and is called the equilibrium measure, cf.\ \cite{Frostman_1935}. When \(F\neq 0\) but \(V=0\) we refer to \cite{Canizo_Carrillo_Patacchini_2015} for existence results for minimisers of $\H$. Even in the case where the energy $\H$ in \eqref{eq:continuum_energy_2} is convex, it may have local minimisers in the Wasserstein \(W_\infty\) or \(W_2\) topology that are not globally minimising, see e.g.\ \cite{Santambrogio_2015} for the definition of the Wasserstein distances. Furthermore, since the discrete energy is never convex (it is defined on a non-convex space) it can have local minimisers or even non-minimising critical points. The study of the mean-field limits of these non-minimising critical points is the object of this article.
We point out that the regularity of critical points of \(\H\) is also of interest. Regularity properties of local minimisers, for the \(W_\infty\) Wasserstein distance, of interaction energies of the form \eqref{eq:continuum_energy_2} were studied in \cite{Balague_Carrillo_Laurent_Raoul_2013,Carrillo_Figalli_Patacchini_2017} (for \(V=0\) and \(g\) less singular than the Coulombian interaction) and  in \cite{Carrillo_Delgadino_Mellet_2016} (for \(V=0\) and Coulombian \(g\)).

 In 2D, mean-field limits of critical points of \(\H_N\) when all the charges \(d_i\) are positive can be described with the help of results by Schochet \cite{Schochet_1996} obtained in the context of point-vortices. Using that result we obtain the convergence of
 stationary solutions of point-vortices to stationary solutions of the Euler equations in the vorticity form. The key argument in Schochet's work requires that all the charges are positive; it relies on  an adaptation of Delort's technique in \cite{Delort_1991} used to pass to the limit in the 2D Euler equations under an assumption of the sign of the vorticity.  The main feature of the present work is to treat the case of sign-changing topological charges \(d_i\). Di Perna-Majda showed in \cite{Diperna_Majda_1988} that it is possible to pass to the limit in the  2D  \textit{stationary} Euler equations without any sign condition on the vorticity measure, cf.\ also \cite{Evans_1990,Bethuel_Ghidaglia_1994, Schochet_1995}.  To the best of our knowledge, characterization of mean-field limits of critical points of \(\H_N\) in the Coulombian case with sign-varying topological charges has not been established yet and our main result can be viewed as an analogue of Di Perna-Majda's result for systems of repulsive-attractive point-vortex systems in 2D. Recently, the existence of weak solutions and the mean-field limit for the (time-dependent) 2D fractional Euler-alignment system with Coulombian singularity was obtained in \cite{F-P_2023}. Remarkably, the issue of diagonal defect measures plays a role here, even though it appears from somewhat different arguments.

\subsection{Statement of the results}\label{sec:results}

To state our results we introduce basic notation. Let $\delta_{ij}$ be the Kronecker delta. For a measurable vector field \(X:\R^d\rightarrow \R^d\) we define the tensor
\begin{equation}\label{eq:stress-grad}
    \lh X,X \rh := \left(2X\otimes X -|X|^2 \text{Id}\right)_{ij} = 2 X_iX_j-|X|^2 \delta_{ij}, \quad \text{ for } i,j=1,\dots,d.
\end{equation}
For any functions \(h,f\in H^1_{\text{loc}}(\R^d)\) we have $\lh\nabla h, \nabla f \rh = [h,f]$, where \( [ h,f]\) is the stress-energy tensor given by  a \(d\times d\) matrix field with entries
\begin{equation}\label{def:stress_tensor}
([h,f])_{ij}:= \p_ih\p_jf+\p_jh\p_if -\nabla h\cdot \nabla f \delta_{ij}  \quad \text{ for } i,j=1,\dots,d.
\end{equation}
Throughout the paper, whenever we can ensure that the vector field $X$ in \eqref{eq:stress-grad} comes from a gradient of some function, say \ $X=\nabla h$, we usually opt to use the notation of the tensor $[h,h]$ in \eqref{def:stress_tensor}. However, at times, to avoid this issue altogether, we simply use \eqref{eq:stress-grad}.

It is easy to check, thanks to the Leibniz rule, that if \(h,f\) are \(\C^2(\R^d)\), then \(\dive ([ h,f])=\Delta h \nabla f+\Delta f\nabla h\). In particular \(\dive [ h,h]=2\Delta h\nabla h\) if \(h\) is \(\C^2(\R^d)\).  This suggests to use \( [ h,h]\) for a weak formulation of \eqref{eq:formal_equilibrium_1}. However \(h_N\) is not in \(H^1_{\text{loc}}(\R^d)\) and thus \([h_N,h_N] \) does not define a distribution. Thus, the rigorous definition of \(\dive \lh \nabla h_N, \nabla h_N \rh\) is unclear. We use the notion of divergence-free in finite part vector fields that we adapt from Sandier-Serfaty \cite[Chapter 7]{Sandier_Serfaty_2007}.
\begin{definition}[Divergence-free in finite part vector fields]\label{def:div_free_fp}
Let \(X\) be a measurable vector field in \(\R^d\), we say that \(X\) is divergence-free in finite part if there exists a family of sets \( (E_\delta)_{\delta>0}\) such that
\begin{enumerate}
\item for any compact set \(K\subset \R^d\), we have \(\lim_{\delta\to 0} \H^1_\infty(K\cap E_\delta)=0,\)
\item for every \(\delta>0,\) the vector field \(X\) is in \(L^1(\R^d\setminus E_\delta)\),
\item for every \(\zeta\) in \(\C^\infty_c(\R^d)\),
\begin{equation*}
\int_{\R^d \setminus \zeta^{-1}(\zeta(E_\delta))} X\cdot \nabla \zeta\dd x=0.
\end{equation*}
\end{enumerate}
\end{definition}
Here \( \H^1_\infty\) denotes the 1-Hausdorff content defined in Definition  \ref{def:Hausdordd_content} below. We say that a tensor is divergence-free in finite part if each of its rows is divergence-free in finite part. Note that this definition mimics in some sense the definition in the sense of distributions of \({\rm vp }\left( \frac{1}{x}\right)\) and \(\rm fp \left(\frac{1}{x^2} \right)\).

Only the criticality of the Coulombian part of the energy necessarily needs to be characterised using the electric potential and the divergence-free in finite part condition. The regular part can be dealt with using more basic techniques. Thus we will treat it separately. Namely, we define the {\it regular part} of the system as
 \begin{equation}\label{def:R_N}
(R_N)_{ij}=  -2 \p_jg \ast \left[\mu_N(\p_i F\ast \mu_N+\p_i V) \right]  \quad \text{ for } i,j=1,\dots, d.
\end{equation}
Then, 
$$\dive R_N =2c_d\mu_N(\nabla F\ast \mu_N+\nabla V).$$ 
We also define the Radon measure \((\nabla F)_\mu\) acting on \(\varphi \in \C^0_c(\R^d)\) by
\begin{equation}\label{eq:false_conv}
    \langle (\nabla F)_\mu, \varphi\rangle :=\frac12 \iint_{\R^d \times \R^d} \nabla F(x-y)\cdot (\varphi(x)-\varphi(y))\dd \mu(x) \dd \mu(y).
\end{equation}
We can check that \((\nabla F)_\mu=\mu (\nabla F \ast \mu)\) whenever the convolution \(\nabla F \ast \mu\) is well-defined and continuous.
Thus, we arrive at three variants of criticality condition for the continuum energy $\H$.
\begin{definition}[Criticality conditions]\label{def:critcond}
    We say that a finite signed Radon measure $\mu$ is a critical point of $\H$, in the sense of 
    \begin{enumerate}
        \item {\bf Weak velocity formulation} if
        $ \lh \nabla h, \nabla h\rh + R $ is divergence-free in finite part, where $R\in L^p_{{\rm loc}}(\R^d;\R^{d\times d})$ for some $1\leq p<\frac{d}{d-1}$ and
        $$\dive R=2c_d((\nabla F)_\mu+\mu \nabla V) \text{ in } \D'(\R^d;\R^d); $$
        \item {\bf Velocity formulation} if the above condition holds, if \(\nabla g \ast \mu \in L^2_{\rm loc}(\R^d;\R^d)\) and 
         $$\dive (\lh \nabla g*\mu,\nabla g*\mu\rh + R)= 0 \text{ in } \D'(\R^d;\R^d);$$
        \item {\bf Vorticity formulation} if for all $\varphi\in \C^\infty_c(\R^d;\R^d)$
        \begin{multline*}
    \frac12 \iint_{(\R^d \times \R^d) \setminus \{x=y\}}\left(\nabla g(x-y)+\nabla F(x-y)\right)\cdot   \left(\varphi(x)-\varphi(y)\right)\dd \mu(x) \dd \mu(y) \\
    + \int_{\R^d}\nabla V(x)\cdot \varphi(x) \dd \mu(x)=0.
\end{multline*}

    \end{enumerate}
\end{definition}

Our first main result concerns the limits of critical points \eqref{eq:3crit} and their criticality in the sense of weak vorticity formulation, cf.\ item (1) in Definition \ref{def:critcond}.

\begin{theorem}\label{th:main_1}
Let \(\vx_N = (x_1,\dots,x_N)\in (\R^d)^N\) be a system of points satisfying \eqref{eq:critical_points}-\eqref{eq:Boundedness_Hamiltonian}.

{\bf Part A} {\rm Convergence of empirical measures and electric potentials}

\noindent
 Then there exist \(\mu\in \M(\R^d)\),  \(R\in L^p_{{\rm loc}}(\R^d;\R^{d\times d})\) for any \(1\leq p< d/(d-1)\) such that, up to a subsequence, 
\begin{align*}
    \mu_N & \to \mu \text{ narrowly} \\
    \nabla h_N &\rightarrow \nabla g \ast \mu \text{ in } L^p_{{\rm loc}}(\R^d;\R^d) \text{ and } R_N \rightarrow R \text{ in } L^p_{{\rm loc}}(\R^d;\R^{d\times d}) \text{ for } 1\leq p< d/(d-1),
\end{align*}
\begin{equation}
    \dive R=2c_d ((\nabla F)_\mu+\mu\nabla V) \text{ in } \D'(\R^d;\R^d).
\end{equation}
Here $\mu_N$, $h_N$ and $R_N$ are the associated empirical measures, electric potentials and non-singular reminders, respectively cf. \eqref{eq:3crit}, \eqref{def:R_N} and \( (\nabla F)_\mu\) is defined in \eqref{eq:false_conv}.

{\bf Part B} {\rm Criticality of \(\mu\) in the weak velocity formulation}

\noindent
    Let $d=2$. Under the assumptions of Theorem \ref{th:main_1} the tensor \(\lh \nabla g*\mu, \nabla g*\mu\rh+R\) is divergence-free in finite part or equivalently, \(\mu\) is a critical point of \(\H\) in the weak velocity formulation (cf.\ Definition \ref{def:div_free_fp}).
\end{theorem}

\noindent
The following result provides conditions ensuring that the critical points established in the above theorem satisfy in fact the augmented criticality condition in item (2) of Definition \ref{def:critcond} and it falls in line with the notion of criticality as a limiting variant of condition \eqref{eq:formal_equilibrium_1}.

\begin{theorem}[Criticality of \(\mu\) in the velocity formulation ]\label{th:main1.5}
    Let $d=2$. Under the assumptions of Theorem \ref{th:main_1}, suppose further that \(\mu\in H^{-1}_{{\rm loc}}(\R^2)\). Then, \(\mu\) is a critical point of \(\H\) in the velocity formulation sense, i.e.
\begin{equation}\label{eq:1.16}
    \dive (\lh \nabla g*\mu,\nabla g*\mu\rh + R)= 0 \text{ in } \D'(\R^2;\R^2). 
\end{equation}
 If additionally \(F\) or \(\mu\) are such that \( \nabla F \ast \mu\) is well-defined and in \(\C^0(\R^2;\R^2)\) and if \( \dive (\nabla F\ast \mu)\) and  \(\Delta V\) belong to \(L^p_{{\rm loc}}(\R^d)\) for some \(p>2\) then \(Y:= \nabla g \ast \mu + \nabla F \ast \mu + \nabla V\) satisfies

\begin{equation}\label{eq:possible_reg}
    \dive \lh Y,Y\rh= 2 (\dive  (\nabla F \ast \mu)+\Delta V)Y.
\end{equation}

\noindent
If  furthermore \(\dive (\nabla F\ast \mu)\), \(\Delta V\) belong to \(L^\infty_{{\rm loc}}(\R^2)\) then \(\nabla g \ast \mu \in L^\infty_{{\rm loc}}(\R^2)\). Lastly, if \(\mu \in L^1(\R^2)\) then \(\mu (\nabla g \ast \mu+\nabla F\ast \mu+\nabla V)=0\) in \(\D'(\R^2;\R^2)\) and actually 

\begin{equation}\label{eq:description}
    \mu= -\left( \dive  (\nabla F \ast \mu )+\Delta V \right) \textbf{1}_{\{|\nabla F \ast \mu+\nabla V +\nabla g \ast \mu|=0\}}=0. 
\end{equation}
\end{theorem}

\begin{remark}[Link with the Euler equations] \rm
As stated previously when \(d=2\), \(F=0\) and \(V=0\), Equation \eqref{eq:critical_points} describes a point-vortex stationary system (without forces). We observe that in that case, if \(\mu\in H^{-1}_{\rm loc}(\R^2)\) the criticality condition obtained in Theorem \ref{th:main1.5} means that \(v:=\nabla^\perp g\ast \mu\in L^2_{\rm{loc}}(\R^2;\R^2)\) is a stationary solution of the Euler equations in fluid mechanics. Indeed, we have  \(\dive (v)=\dive(\nabla^\perp g\ast \mu)=0\) and  
\begin{align*}
    2 v\otimes v -|v|^2 \rm{Id}&= [\nabla^\perp g \ast \mu, \nabla^\perp g \ast \mu]=-[\nabla g \ast \mu,\nabla g \ast \mu].
\end{align*}
Hence we arrive at \(\dive (v\otimes v)=-\nabla p\) in \(\D'(\R^2;\R^2)\) with \(p=\frac{|v|^2}{2}\) and \(v\) satisfies the Euler equation with the pressure given by \(|v|^2/2\).

\end{remark}
\begin{remark}[Growth conditions on $F$]\label{rem:otherass}\rm
    We point out that \(\nabla F \ast \mu\) is well-defined and in \(\C^0\) when either \(\mu\) has compact support or \(\nabla F \) is bounded. The condition that \(\dive (\nabla F \ast \mu) \in L^p_{\rm loc}(\R^d)\) can be achieved for example when \(\Delta F\in L^p(\R^d) \). More generally, if $|\nabla F(x)|\lesssim |x|^p +1$, then for \(\nabla F \ast \mu\) to be well-defined, one requires uniform integrability of $p$-th moments of $\mu_N$.
\end{remark}

\noindent
Our third main result employs the vorticity formulation of \eqref{eq:critical_points} in the sense of item (3) in Definition \ref{def:critcond} inspired by the one used in \cite{Delort_1991,Schochet_1995} for the 2D Euler equations.

\begin{theorem}[Vorticity criticality condition]\label{th:main_2}
    Let \(\vx_N= (x_1,\dots,x_N)\in (\R^2)^N\) and $\mu_N$ (cf. \eqref{eq:3crit}) be a sequence satisfying \eqref{eq:critical_points}. Then, up to a subsequence, $\mu_N$ converges narrowly to some \(\mu\in \M_f(\R^2)\) and if \(\mu\) is  non-atomic then
    \begin{equation}\label{eq:limit_vorticity_form}
        \frac12 \iint_{\R^2\times \R^2}  \left(\nabla g(x-y)+\nabla F(x-y)\right) \cdot (\varphi(x)-\varphi(y)) \dd \mu(x) \dd \mu(y) +\int_{\R^2} \nabla V(x)\cdot \varphi(x) \dd \mu(x)=0
    \end{equation}
    for all \(\varphi \in \C^\infty_c(\R^2;\R^2)\).
\end{theorem}

\begin{remark}[Equivalence between velocity and vorticity formulations]\rm
 Any finite signed Radon measure \(\mu \in H^{-1}_{\text{loc}}(\R^2)\)  is non-atomic, cf.\ e.g.\ \cite[Lemma 1.2.5]{Delort_1991}. Then, the criticality condition obtained in Theorem \ref{th:main1.5} and the one obtained in Theorem \ref{th:main_2} are equivalent, cf.\ Proposition \ref{prop:equivalence} below. However the criticality condition in Theorems \ref{th:main_1}  is valid even if \(\mu\) is merely a finite sum of Dirac masses. In such a case it states that the flux of \(\lh\nabla h,\nabla h \rh\) is equal to \(-2\nabla F \ast \mu (x_i) -2\nabla V(x_i) \) around the point \(x_i\), cf.\ Proposition \ref{prop:equiv_div_fp_zero_flux}. 
\end{remark}
 The advantage of the formulation employing electric potentials, e.g.\ \eqref{eq:possible_reg}, is that it can be used to study the regularity of the limiting measures. Indeed it was shown in \cite{Rodiac_2016} that measures satisfying \eqref{eq:possible_reg} with \(F=V=0\) are locally supported by curves which are zero sets of a harmonic function (possibly multi-valued). We also refer to \cite{Rodiac_2019} for the study of regularity of measures satisfying conditions very similar to \eqref{eq:possible_reg} in the context of Ginzburg-Landau vortices. 

Our last main result deals with the limit of stable critical points of \(\mathcal{H}_N\). This time we assume that all the charges are non-negative and we can show that, under that assumption, stability passes to the limit in the following sense. 
\begin{definition}[Stable critical points]\label{def:stability}
Assume that \(F, V\in \C^2(\R^2)\). We say that \( (x_1,\dots,x_N)\) is a stable critical point of \(\mathcal{H}_N\) if it satisfies \eqref{eq:critical_points} and if 
\begin{align}
    \sum_{i\neq j=1}^N d_i d_j ( D^2g(x_i-x_j)+D^2F(x_i-x_j))+ d_i D^2V(x_i) \geq 0 &\text{ for all } i\in\{1,\dots, N\}, \label{eq:stability_1}\\
    d_i d_j D^2g(x_i-x_j)+D^2F(x_i-x_j) \geq 0 & \text{ for all } i\neq j\in\{1,\dots, N\}. \label{eq:stability_2}
\end{align}
\end{definition}
For the sake of the following theorem, we recall that a symmetric matrix $A$ is positive-definite, denoted \(A\geq 0\), provided that \(x^TAx\geq 0\) for all \(x\in \R^d\).
\begin{theorem}\label{th:main3}
    
    Assume that \(F,V\) are in \(\C^2(\R^2)\). Let \(\vx_N= (x_1,\dots,x_N)\in (\R^2)^N\) and $\mu_N$ be as in \eqref{eq:3crit} and satisfy \eqref{eq:critical_points}-\eqref{eq:Boundedness_Hamiltonian}. Assume that \(d_i\geq 0\) for all \(i=1,\dots,N\). Then, up to a subsequence, \(\mu_N \to \mu\) narrowly,  and if \(\mu\) is non-atomic then
    \begin{equation}\label{eq:stability_limit}
        \iint_{\R^2 \times \R^2} (\varphi(x)-\varphi(y))^T\left(D^2g(x-y)+ D^2F(x-y)\right) (\varphi(x)-\varphi(y)) \dd \mu(x) \dd \mu(y) \geq 0
    \end{equation}
    for all \(\varphi \in \C^\infty_c(\R^2;\R^2).\)
\end{theorem}

Condition \eqref{eq:stability_limit} can be interpreted as follows: when all the charges are positive, then stable critical points of $\H_N$ converge to stable critical points of $\H$, under the assumption that $V=0$ (which, remarkably, does not break the assumption \eqref{A2} in dimension 2). Passing to the limit in the stability condition for Ginzburg-Landau vortices was recently investigated in \cite{Rodiac_2024} under an assumption of convergence of energies.

\subsection{Sketch  of the proofs and organization of the paper}

The paper is organized as follows. Section \ref{sec:prelim} is devoted to the introduction of the preliminary information. First, in Section \ref{sec:reformulation}, we focus on the reformulation of equation \eqref{eq:critical_points} in terms of the electric stress tensor $[h_N,h_N]$, cf. \eqref{eq:3crit} and \eqref{def:stress_tensor}.
The reformulation involves a vanishing  flux condition for \(\lh\nabla h_N, \nabla h_N\rh+R_N\) which can be reformulated in terms of the divergence-free in finite part notion. We prove it and collect some properties related to the divergence-free in finite part condition in Section \ref{sec:divergence_free_finite_part}. In Section \ref{sec:convergence} we prove various types of compactness for the sequences of critical points in terms of empirical measures and associated electric potentials \eqref{eq:3crit}.
Particularly, we show that 
since, up to a subsequence, \(\mu_N\) converges towards \(\mu\) narrowly, thanks to the Sobolev embedding \(W^{1,q} \hookrightarrow \C^0\) for \(q>d\), we infer that, up to a subsequence, \(\mu_N \rightarrow \mu\) in \(W^{-1,p}(\R^d)\) for any \(1\leq p<d/(d-1)\). Elliptic regularity theory, or better said, the study of the convolution with \(\nabla g\) shows that \(\nabla g \ast \mu_N \rightarrow \nabla g \ast \mu\) in \(W^{1,p}_{\text{loc}}(\R^d,\R^d)\) for \(1\leq p<d/(d-1)\).
In Section \ref{sec:velocity_form}, key in our article, we refine this convergence and show that  \(\nabla g \ast \mu_N\) converges, up to a subsequence, towards \(\nabla g \ast \mu\) strongly in \(L^2\) except on some sets with arbitrary small \(1\)-capacity, cf.\ Proposition \ref{prop:conv_capacity}. This fact allows us to pass to the limit in the divergence-free in finite parts condition and to obtain part B of Theorem \ref{th:main_1}. In order to pass to the limit  we use  an argument involving the coarea formula as it is done in \cite[Chapter 13]{Sandier_Serfaty_2007}, see also \cite{Sandier_Serfaty_2003} and \cite{Bethuel_Ghidaglia_1994}. We augment the criticality conditions under higher regularity of the limit and prove Theorem \ref{th:main1.5} in Section \ref{sec:aug}.
Then, in Section \ref{sec:vorticity_form}, we  pass to the limit in the vorticity form of \eqref{eq:critical_points} proving Theorem \ref{th:main_2}. We first show that \(\mu_N\) satisfies an equation of the form \eqref{eq:vorticity_form}. Then since the kernel against which we integrate \(\mu_N \otimes \mu_N\) is bounded, employing a defect measure argument, we show that the limiting measure satisfies an equation involving a defect measure carried by the diagonal. A suitable choice of test functions $\varphi$ reveals that this defect measure actually vanishes. The arguments in this section follow the ones in \cite{Schochet_1995}.
In the final Section \ref{sec:final}, we address our last main result, again using arguments  related to the vorticity form, cf. Section \ref{sec:vorticity_form}. The kernel appearing in the stability condition for \(\vx_N\) is bounded. However, we cannot use the same argument of defect measure as before since we have to pass to the limit in an inequality and not an equality. What allows us to succeed is that, in the case where the empirical measures are non-negative and if we integrate the product of these measures against a bounded kernel (continuous outside the diagonal), we can show directly that the desired quantity passes to the limit. This argument is originally due to Delort \cite{Delort_1991}.

\subsection{Notation}

By \(\mathcal{M}_f(\R^d)\) we denote the set of all finite signed Radon measures on \(\R^d\), 
and by \(\mathcal{M}_f^+(\R^d)\) -- its subset consisting of non-negative measures.
 We say that a family ${\mathcal F}\subset \mathcal{M}_f(\R^d)$ is {\it bounded} if
 $$ \sup_{\mu\in{\mathcal F}} |\mu|_{TV} < +\infty. $$
    where \(|\mu|_{TV}:=\sup \{\int_{\R^d} \varphi \dd \mu; \varphi \in \C^0_c(\R^d), \|\varphi\|_{\infty}\leq 1\}\) is the total variation of \(\mu\).  We say that it is {\it tight} if for all $\varepsilon>0$ there exists a compact set $K\subset\R^d$, such that
    $$ \sup_{\mu\in{\mathcal F}}|\mu|(\R^d\setminus K)<\varepsilon. $$
    Here $|\mu|$ denotes the variation of measure $\mu$, more precisely if \(\mu=\mu^+-\mu^{-}\) with \(\mu^+,\mu^{-}\in \mathcal{M}^+(\R^d)\) then \(|\mu|=\mu^++\mu^-\). Finally, we say that a sequence of measures \((\mu_N)_N \in \mathcal{M}_f(\R^d)\) converges {\it narrowly} to some \(\mu \in \mathcal{M}_f(\R^d)\) if 
    \begin{equation}
    \int_{\R^d} \varphi \dd \mu_N \rightarrow \int_{\R^d} \varphi \dd \mu \quad \forall \varphi \in \C^0_b(\R^d),
\end{equation}
where \(\C^0_b(\R^d)\) denotes the set of continuous bounded functions on \(\R^d\). 
Note that the narrow convergence is stronger than the weak-* (or vague) convergence tested by continuous compactly supported (or vanishing at infinity) functions. We also observe that the weak-* and vague convergence are equivalent in \(\R^d\).
\noindent
Prokhorov theorem for signed Radon measures, see e.g.\ \cite[Theorem 1.4.11]{Bogachev_2018}, states that a sequence in $\M_f(\R^d)$ is tight if and only if it has a narrowly convergent subsequence.

We denote the product measure of two measures \(\mu, \nu \in \mathcal{M}_f(\R^d)\) as \(\mu \otimes \nu\), which belongs to the space $\M_f(\R^{2d})$. 
We denote the inner product of two vectors by \(u\cdot v\). For two matrices their Frobenius inner product is defined by \(A:B=\tr (A^TB)\).

For simplicity of notation, we sometimes write $A\lesssim B $
to indicate that there exists a positive ``harmless'' constant $C$, such that $A\leq CB$.

\subsection{Acknowledgments:} The research of R.R is part of the project \\
No. 2021/43/P/ST1/01501 co-funded by the National Science Centre and the European Union Framework Programme for Research and Innovation Horizon 2020 under the Marie Skłodowska-Curie grant agreement No. 945339. For the purpose of Open Access, the author has applied a CC-BY public copyright licence to any Author Accepted Manuscript (AAM) version arising from this submission. J.P. was partially supported by the Polish National Science Centre’s Grant No. 2018/31/D/ST1/02313 (SONATA) and by University of Warsaw program IDUB Nowe Idee 3A.

\section{Preliminaries}\label{sec:prelim}
We dedicate the preliminaries to two main issues. First,  we reformulate the criticality condition  \eqref{eq:critical_points} in terms of the electric potential $h_N$, and second, we provide a general overview of the notion of divergence-free in finite part vector fields.

\subsection{Reformulation of the criticality condition for the particle system}\label{sec:reformulation}

In this section, we reformulate the system \eqref{eq:critical_points} into  a system of \(N\) flux conditions for the electric potential $h_N$ in \eqref{eq:3crit}. This is inspired by the characterization of criticality for the renormalized energy in the Ginzburg-Landau model due to Bethuel-Brezis-H\'elein \cite{Bethuel_Brezis_Helein_1994}. 

Recall the tensor $\lh \nabla h,\nabla f\rh = [h,f]$ in \eqref{eq:stress-grad} and \eqref{def:stress_tensor}, and that throughout the paper $g$ is the fundamental solution of the Laplace equation. We start with the following two technical lemmas.
\begin{lemma}\label{lem:chain_rule_stress_tensor}
Let \(\O\) be an open set in \(\R^d\), for \(d\geq 2\) and fix \(1\leq p \leq  +\infty\) with \(p'\) such that \(\frac{1}{p}+\frac{1}{p'}=1\). Then 
\begin{equation}\label{eq:stress_energy_equality}
\dive [h,f]= \Delta h\nabla f+\Delta f\nabla h,
\end{equation}
provided that any of the following conditions is satisfied:
\begin{enumerate}
    \item the functions \(h,f\) belong to \(\C^2(\O)\),
    \item we have \(\nabla h \in L^p(\O;\R^d), \nabla f\in L^{p'}(\O;\R^d)\), \(\Delta h \in L^p(\O), \Delta f\in L^{p'}(\O)\),
    \item we have \(\nabla h \in L^p(\O;\R^d), \nabla f\in \C^0(\O;\R^d)\) and \(\Delta h \in \mathcal{M}(\O), \Delta f \in L^{p'}(\O)\).
\end{enumerate}

\end{lemma}

\begin{proof}
    If \(h,f\) are in \(\C^2(\O)\), identity \eqref{eq:stress_energy_equality} is proven thanks to the chain rule. For \(i=1,\dots,d\), we have
    \begin{align*}
        \left(\dive [h,f] \right)_i &= \sum_{j=1}^d \p_j \left(\p_ih\p_jf+\p_jh\p_if -\nabla h\cdot \nabla f \delta_{ij} \right) \\
        &= \sum_{j=1}^d \bigl( \p^2_{ij} h\p_jf+\p^2_{jj}h\p_i f+\p_ih \p^2_{jj}f+\p_jh \p^2_{ij} f \bigr) \\
        &\quad \quad -\sum_{j=1}^d\sum_{k=1}^d \left(\p^2_{jk} h\p_kf+\p_kh\p^2_{jk} f\right)\delta_{ij} \\
        &=\Delta h \nabla f +\Delta f \nabla h.
    \end{align*}
    Now, in the case of item (2) the result is proven by approximation. Indeed, proceeding as in \cite[Theorem 4.2]{Evans_Gariepy_2015}, we can find \(h_k,f_k\in \C^\infty(\O)\) such that 
    \begin{align*}
        \nabla h_k \rightarrow \nabla h  \text{ in } L^p(\O;\R^d), & \quad  \nabla f_k \rightarrow \nabla f \text{ in }  L^{p'}(\O;\R^d) \\
        \Delta h_k \rightarrow \Delta h \text{ in } L^p(\O), & \quad \Delta h_k \rightarrow \Delta f \text{ in } L^{p'}(\O).
    \end{align*}
    Hölder's inequality shows that \( [h_k,f_k] \rightarrow [h,f]\) in \(L^1(\O;\R^{d\times d})\), and thus \(\dive [h_k,f_k] \rightharpoonup \dive [h,f]\) in \(\D'(\O;\R^d)\). But the same Hölder's inequality shows that \(\Delta h_k \nabla f_k+\Delta f_k \nabla h_k \rightarrow \Delta h\nabla f+\Delta f\nabla h\) in \(L^1(\O;\R^d)\). Hence equality \eqref{eq:stress_energy_equality} holds. The case when \(\Delta h \in \mathcal{M}(\O), \nabla f\in \C^0(\O)\) is also proven by approximation, since then \(\nabla f_k \rightarrow \nabla f\) locally uniformly in \(\O\) and \(\Delta h_k \to \Delta h\) weakly-* as signed Radon measures. This is sufficient to prove that \(\Delta h_k \nabla f_k \rightharpoonup \Delta h \nabla f\) in \(\D'(\O;\R^d)\).
    \end{proof}

\begin{lemma}\label{lem:cancelation}
 Let \(d\geq 2\), and \(x_0\in\R^d\). For any \(\delta>0\), if \(\nu\) denotes the outward unit normal to  \(B_\delta(x_0)\) then
 \begin{equation*}
 \int_{\p B_\delta (x_0)} [ g(\cdot-x_0),g(\cdot-x_0)]\nu \dd \mathcal{H}^{d-1}=0.
 \end{equation*}
\end{lemma}

\begin{proof}
For \(k=1,\dots,d\) and for \(x\neq x_0\) we compute that 
\begin{equation*}
\p_kg(x-x_0)= \begin{cases}
-\frac{x^k-x_0^k}{|x-x_0|^{2}} &\text{ if } d=2 \\
    -(d-2)\frac{x^k-x_0^k}{|x-x_0|^{d}} & \text{ if } d\geq 3,
\end{cases}
 \quad  |\nabla g(\cdot-x_0)|^2= \frac{1}{|x-x_0|^{2(d-1)}} \ \text{for all } d\geq 2.
\end{equation*}
The outward unit normal to \(\p B_\delta(x_0)\) is given by \(\nu(x)=\frac{x-x_0}{|x-x_0|}\). For \(d\geq 3\), the \(k\)th-component of the vector \[\int_{\p B_\delta (x_0)} [ g(\cdot-x_0),g(\cdot-x_0)]\nu \dd \mathcal{H}^{d-1}\] is given by 
\begin{align*}
\int_{\p B_\delta(x_0)} 2(d-2)^2 \sum_{\ell=1}^d\frac{(x^k-x_0^k)(x^\ell-x_0^\ell)}{|x-x_0|^{2d}}\frac{x^\ell-x_0^\ell}{|x-x_0|}-\frac{(d-2)^2}{|x-x_0|^{2(d-1)}}\frac{x^k-x_0^k}{|x-x_0|} \dd \H^{d-1}\\
=\int_{\p B_\delta(x_0)} \frac{(d-2)^2}{|x-x_0|^{2(d-1)}}\frac{x^k-x_0^k}{|x-x_0|} \dd \H^{d-1} \\
=\frac{(d-2)^2}{\delta ^{2(d-1)}}\int_{\p B_\delta(x_0)} \nu^k \dd \H^{d-1}=0.
\end{align*}
The computation for \(d=2\) is similar and can be found in \cite[Example 2.8]{Rodiac_2016}. 
\end{proof}

We are now ready to state the first reformulation of \eqref{eq:critical_points} in terms of the electric potential \(h_N\).

\begin{proposition}\label{prop:flux_condition}
Let \( \vx_N = (x_1,...,x_N)\in (\R^d)^N\), and \(\mu_N\), \(h_N\) be the associated empirical measure  and the electric potential, respectively, cf.\ \eqref{eq:3crit}. Finally, let \(R_N\) be regular reminder defined in \eqref{def:R_N}. Then 
\begin{equation}\label{eq:divergence_free_condition_1st}
\dive \left( [h_N,h_N]+R_N\right)=0 \text{ in } \D'(\R^d\setminus \{x_1,\dots,x_N\};\R^d).
\end{equation}
Moreover \( \vx_N\) satisfies \eqref{eq:critical_points} if and only if 
\begin{equation}\label{eq:flux_condition}
\int_{\p B_\delta (x_i)} \left(  [h_N,h_N]+R_N \right)\nu \dd \mathcal{H}^{d-1}=0, \quad \forall 0<\delta< \min_{i\neq j} |x_i-x_j|,\ \forall i=1,\dots,N.
\end{equation}
\end{proposition}

\begin{proof}
We only show the \(d \geq 3\) case ($d=2$ is similar).
We start by observing that the vector-valued function 
$$-2g \ast \left[\mu_N(\p_k F\ast \mu_N+\p_k V) \right]$$
is in \(\C^\infty(\R^d\setminus \{x_1,\dots,x_N\};\R^d)\) since it is harmonic outside \(\{x_1,\dots,x_N\}\) (it is a convolution with the fundamental solution of the Laplace equation). Hence its derivative \(R_N\) is also in \(\C^\infty(\R^d\setminus \{x_1,\dots,x_N\};\R^{d\times d})\) and, by definition, \(\dive R_N=2c_d\mu_N(\nabla F\ast \mu_N+\nabla V)=0\) in \(\D'(\R^d\setminus \{x_1,\dots,x_N\};\R^d)\). Similarly, \(h_N\) is harmonic in \(\R^d\setminus\{x_1,\dots,x_N\}\), and hence, from Lemma \ref{lem:chain_rule_stress_tensor}, we see that \(\dive [h_N,h_N]=2\Delta h_N \nabla h_N =0\) in \(\D'(\R^d\setminus \{x_1,\dots,x_N\};\R^d)\). Thus the proof of \eqref{eq:divergence_free_condition_1st} is concluded.

In order to prove the second assertion \eqref{eq:flux_condition} we aim to show that 
\begin{equation}\label{eq:flux_condition2}
\begin{split}
    \int_{\p B_\delta(x_i)} &([h_N,h_N]+R_N)\nu \dd \H^{d-1}\\
    &=\frac{2c_d}{M_N^2}\sum_{i\neq j} d_i d_j\left[ \nabla g(x_i-x_j)+F(x_i-x_j) \right] +\frac{2c_d}{M_N}d_i\nabla V(x_i)
\end{split}
\end{equation}
for all $\delta$ as in \eqref{eq:flux_condition}. Then the right-hand side above vanishes if and only if $\vx_N$ satisfies \eqref{eq:critical_points}.

\noindent
Denote \(q_N:=\inf_{i\neq j} |x_i-x_j|\). Since \(R_N\in \C^\infty(\R^d\setminus \{x_1,\dots,x_N\};\R^{d\times d})\) and \(\dive R_N \in \mathcal{M}(\R^d;\R^d)\), the divergence theorem implies that, for \(\delta<q_N\), we have
\begin{align}
    \int_{\p B_{\delta}(x_i)} R_N \nu \dd \H^{d-1} &= \int_{B_\delta (x_i)} \dd (\dive R_N) = \int_{B_\delta (x_i)} 2(\nabla F\ast \mu_N +\nabla V) \dd \mu_N \nonumber\\
    & = \frac{2}{M_N^2} \sum_{j\neq i} d_i d_j F(x_i-x_j) +\frac{2}{M_N} d_i\nabla V(x_i). \label{eq:first_part}
\end{align}

\noindent
The function
\begin{equation}\label{H_NI}
H_{N,i}:= \frac{1}{M_N}\sum_{j\neq i} d_j g(\cdot - x_j)
\end{equation}
is harmonic, and hence smooth in \( B_\delta(x_i)\) for \(\delta <q_N\). Moreover, since \(h_N=\frac{1}{M_N} d_ig(\cdot-x_i)+H_{N,i}\), we obtain
\begin{equation}\label{eq:decompo_h_N}
    [h_N,h_N] =\frac{d_i^2}{M_N^2}[g(\cdot-x_i),g(\cdot-x_i)]+\frac{2d_i}{M_N}[g(\cdot-x_i),H_{N,i}]+[H_{N,i},H_{N,i}].
\end{equation}
However \( \int_{B_\delta(x_i)}[g(\cdot-x_i),g(\cdot-x_i)]\nu \dd \H^{d-1}=0 \) by Lemma \ref{lem:cancelation}.  We also have, using the divergence theorem and Lemma \ref{lem:chain_rule_stress_tensor},
\begin{align*}
    \int_{\p B_\delta(x_i)} [H_{N,i},H_{N,i}]\nu \dd \H^{d-1}& =\int_{B_\delta(x_i)} \dive [H_{N,i},H_{N,i}] \dd x = \int_{B_\delta(x_i)} 2\Delta H_{N,i}\nabla H_{N,i} \dd x=0.
\end{align*}

Now we use Einstein's summation convention for \(k=1,\dots,d\) to compute
\begin{equation}\label{eq:pre_delta}
\begin{split}
\int_{\p B_\delta(x_i)} &([g(\cdot-x_i),H_{N,i}]\nu)_k  \dd \H^{d-1}\\
&=\int_{\p B_\delta(x_i)} \left( \p_kg(x-x_i)\p_\ell H_{N,i}(x)\nu^\ell(x) +\p_\ell g(x-x_i)\p_kH_i(x)\nu^\ell (x)\right)\dd \H^{d-1}  \\
&\qquad -\int_{\p B_\delta(x_i)} \nabla g(x-x_i)\cdot \nabla H_{N,i}(x)\nu^k(x)\dd \H^{d-1} \\
&= -(d-2)\int_{B_\delta(x_i)} \left(\frac{x^k-x_i^k}{|x-x_i|^d}\frac{x^\ell-x_i^\ell}{|x-x_i|}\p_\ell H_{N,i}(x) - \frac{(x^\ell-x_i^\ell)^2}{|x-x_i|^{d+1}}\p_k H_{N,i}(x)\right) \dd x \\
&\qquad+(d-2)\int_{\p B_\delta(x_i)} \frac{(x^k-x_i^k)(x^\ell-x_i^\ell)}{|x-x_i|^{d+1}}\p_\ell H_{N,i}(x)\dd \H^{d-1}  \\
&=\frac{(d-2)}{\delta^{d-1}} \int_{\p B_\delta(x_i)} \p_kH_{N,i}(x)\dd \H^{d-1} . 
\end{split}
\end{equation}

\noindent
On the other hand, since we already know that \(\dive [h_N,h_N]=0\) in \(\D'(\R^d\setminus \{x_1,\dots,x_N\};\R^d)\), the quantity \(\int_{\p B_\delta(x_i)} ([g(\cdot-x_i),H_{N,i}]\nu)_k\dd \H^{d-1}\) is independent of \(\delta\) for \(0<\delta<q_N\), and it is equal to its limit as \(\delta \to 0\). By the smoothness of \(H_{N,i}\) we find that 
\begin{align}\label{eq:final_cond}
\lim_{\delta \to 0} \frac{(d-2)}{\delta^{d-1}} \int_{\p B_\delta(x_i)} \p_kH_{N,i}(x) \dd \H^{d-1}& = (d-2)\mathcal{H}^{d-1}(\mathbb{S}^{d-1}) \p_kH_{N,i}(x_i) \nonumber  \\
&=(d-2)\mathcal{H}^{d-1}(\mathbb{S}^{d-1}) \frac{1}{M_N}\sum_{j\neq i} d_j\p_k g(x_j-x_i) \nonumber \\
&=\frac{c_d}{M_N}\sum_{j\neq i} d_j\p_k g(x_j-x_i),
\end{align}

\noindent
where $c_d$ is such that $-\Delta g=c_d\delta_0$. This is valid for any \(k=1,\dots,d\) and any \(i=1,\dots,N\). Thus, using \eqref{eq:first_part},\eqref{eq:decompo_h_N} and \eqref{eq:final_cond} we obtain \eqref{eq:flux_condition2} and conclude the proof.
\end{proof}

We end this section by providing a corollary, which condenses all of the above results and applies them directly to the case, which is the most relevant from our point of view, i.e.\ to sequences of critical points of $\H_N$.

\begin{corollary}\label{cor:flux_condition_imp}
   Let \(\vx_N = (x_1,\dots,x_N)\), $\mu_N$ and $h_N$ satisfy \eqref{eq:critical_points} and \eqref{eq:3crit}. Assume further that \(R_N\) is the regular part of the system defined in \eqref{def:R_N}, \(U\) is a \(\C^1\) bounded open set in \(\R^d\) and there exists \(\rho>0\) such that \(B_\rho(x_i)\cap B_\rho(x_j)=\emptyset\) for all \(i\neq j\) and \(\p U \cap \bigcup_{i=1}^N B_\rho (x_i)=\emptyset\). Then
   \begin{equation}
       \int_{\p U} \left( [h_N,h_N]+R_N\right)\nu \dd \mathcal{H}^{d-1} =0.
   \end{equation}
\end{corollary}

\begin{proof}
    Since we know that \(\dive ([h_N,h_N]+R_N)=0\) in \(\D'(U\setminus \bigcup_{i=1}^N B_\rho (x_i);\R^d)\) we can apply the divergence theorem and  Proposition \ref{prop:flux_condition} to obtain  
    \begin{align*}
        \int_{\p U} \left( [h_N,h_N]+R_N\right)\nu \dd \mathcal{H}^{d-1}=- \sum_{\{i: B_\rho(x_i)\in U\}} \int_{  \p B_\rho(x_i)} \left( [h_N,h_N]+R_N\right)\nu \dd \mathcal{H}^{d-1} =0;
    \end{align*}
\end{proof}

\subsection{Divergence-free in finite part vector fields}\label{sec:divergence_free_finite_part}

Here we introduce basic information on divergence-free in finite parts vector fields, cf.\ \cite[Chapter 13]{Sandier_Serfaty_2007}. We show that this notion is equivalent to the distributional divergence-free condition for vector fields in \(L^1_{\text{loc}}\). We also prove that the reformulation of the  criticality condition \eqref{eq:critical_points} in Proposition \ref{prop:flux_condition} can be expressed in terms of the divergence-free in finite parts condition. We recall Definition \ref{def:div_free_fp} and we start by introducing the Hausdorff content.
\begin{definition}\label{def:Hausdordd_content}
Let \(E\) be a set in \(\R^d\) and \(k\in \mathbb{N}\), the \(k\)-dimensional Hausdorff content of \(E\) is defined by 
\begin{equation*}
\H^k_\infty(E)=\inf \left\{ \sum_{j=1}^{ +\infty} r_j^k: E\subset \bigcup_{j=1}^{ +\infty} B_{r_i}(x_i)\right\},
\end{equation*}
where the infimum is taken over all countable collections of balls \(\{B_{r_i}(x_i)\}_{i\in \mathbb{N}}\) that cover \(E\).
\end{definition}

\begin{remark}\rm
 Observe that for any set \(E\subset \R^d\) and any integer \(k\), we have \(\H^k_\infty(E)\leq \H^k(E)\), where \(\H^k\) denotes the \(k\)-dimensional Hausdorff measure. Moreover, if \(\zeta:\R^d\rightarrow \R\) is a Lipschitz function then
\begin{equation}\label{eq:Lipschitz_rule}
    \mathcal{L}^1(\zeta(E))\leq 2\|\zeta\|_{W^{1,\infty}}\H^1_\infty(E) \text{ for any } E\subset \R^d.
\end{equation}
\end{remark}
 The following proposition shows that the above notion generalizes the classical divergence-free condition for \(L^1_{\text{loc}}\) vector fields.
\begin{proposition}\label{prop:div_free_fp+L1_div_free}
Let \(X:\R^d\rightarrow \R^d\) be a vector field that is divergence-free in finite part. We assume that \(X\in L^1_{\loc} (\R^d\setminus E;\R^d)\) for \(E\subset \R^d\), then for any \(\zeta \in \C^\infty_c(\R^d)\),
\begin{equation*}
\int_{\R^d \setminus \zeta^{-1}(\zeta(E))} X\cdot \nabla \zeta \dd x=0.
\end{equation*}
In particular if \(X\in L^1_{\loc}(\R^d;\R^d)\) then we can take \(E=\emptyset\) obtaining \( \dive X=0\) in \(\D'(\R^d)\).
\end{proposition}

\begin{proof}
Let \(\zeta\) be a function in \(\C^\infty_c(\R^d)\), and let \(C:=\{x\in \R^d:\ \nabla \zeta (x)=0\}\). We have
\begin{align*}
\int_{\R^d \setminus \zeta^{-1}(\zeta(E))} X\cdot \nabla \zeta \dd x=\int_{\R^d \setminus (\zeta^{-1}(\zeta(E))\cup C)} X\cdot \nabla \zeta \dd x=\int_{\R^d \setminus(\zeta^{-1}(\zeta(E))\cup C)}  X\cdot \frac{\nabla \zeta}{|\nabla \zeta|}|\nabla \zeta|\dd x.
\end{align*}
Applying the coarea formula, we infer that for almost all  regular values \(t\) of \(\zeta\) belonging to \(\R\setminus \zeta(E)\) the function 
\begin{equation*}
f(t):= \int_{\{x:\zeta(x)=t\}} X\cdot \frac{\nabla \zeta}{|\nabla\zeta|}\dd \H^{d-1}
\end{equation*}
is well-defined and in \(L^1(\R\setminus \zeta(E \cup C))\) and 
\begin{equation}\label{eq:coarea_1n}
\int_{\R^d\setminus\zeta^{-1}(\zeta(E))} X\cdot \nabla \zeta\dd x=\int_{\R\setminus \zeta(E)} \int_{\{x:\zeta(x)=t\}} X\cdot \frac{\nabla \zeta}{|\nabla \zeta|}\dd \H^{d-1}\dd t,
\end{equation}
where we have used Sard's theorem stating that \(\mathcal{L}^1(\zeta(C))=0\). We assume that \(X\) is divergence-free in finite part, hence there exist sets \(E_\delta\) as in Definition \ref{def:div_free_fp}. As our next step, let us prove that 
\begin{equation}\label{eq:frowne0}
f(t)=0\quad \text{ for a.e. } t \in \R \setminus(\zeta (E)\cup \zeta(E_\delta)).
\end{equation}
Indeed, let \(\xi:\R\rightarrow \R\) be smooth, injective with \(\xi(0)=0\). Applying the coarea formula we obtain
\begin{align*}
\int_{\R^d\setminus (\xi\circ \zeta)^{-1}(\xi\circ\zeta(E_\delta))} X\cdot \nabla (\xi\circ \zeta)\dd x &=\int_{\R^d\setminus \zeta^{-1}(\zeta(E_\delta))} X(x)\cdot \nabla \zeta(x) \xi'(\zeta(x)) \dd x\\
&=\int_{\{t: t\notin \zeta(E_\delta)\}} \xi'(t)f(t) \dd t.
\end{align*}

\noindent
Since \(\xi\circ \zeta\in \C^\infty_c(\R^d)\) we can apply the divergence-free in finite part condition of \(X\)  to deduce that 
\[ \int_{\{t: t\notin \zeta(E_\delta)\} } \xi'(t) f(t) \dd t=0 \quad \text{ for any } \xi\in \C^\infty(\R), \ \xi(0)=0,\ \xi \text{ injective}.\]
In particular 
\[ \int_{\{t: t\notin \zeta(E_\delta)\} } \beta(t) f(t) \dd t=0 \quad \text{ for any positive } \beta\in \C^\infty(\R),\]
which can be extended by approximation to nonnegative continuous $\beta$ and further to all continuous functions $\beta$ (by expressing $\beta$ as a sum of its positive and negative parts). Ultimately \eqref{eq:frowne0} is proved.

\noindent
By \eqref{eq:coarea_1n} and  \eqref{eq:frowne0}, we have
\begin{align*}
\int_{\R^d\setminus \zeta^{-1}(\zeta(E))} X\cdot \nabla \zeta\dd x&=\int_{\R\setminus \zeta(E)} f(t) \dd t \\
&=\int_{\R\setminus (\zeta(E)\cup \zeta(E_\delta))} f(t) \dd t+\int_{(\R\setminus \zeta(E))\cap \zeta(E_\delta)} f(t) \dd t\\
&=\int_{(\R\setminus \zeta(E))\cap \zeta(E_\delta)} f(t) \dd t.
\end{align*}
Recall that, according to Definition \ref{def:div_free_fp}, \(\lim_{\delta \to 0} \mathcal{H}^1_\infty(E_\delta \cap \supp \zeta)=0\), and by using \eqref{eq:Lipschitz_rule}, we find that \(\lim_{\delta \to 0} \mathcal{L}^1(\zeta(E_\delta))=0\). Thus, since \(f\in L^1(\R\setminus\zeta(E))\) we have 
$$\lim_{\delta \to 0} \int_{(\R\setminus \zeta(E)) \cap \zeta(E_\delta)} f(t) \dd t=0.$$
We obtain the conclusion \(\int_{\R^d\setminus \zeta^{-1}(\zeta(E))} X \cdot \nabla \zeta\dd x=0.\)

\end{proof}

Next, we prove that the flux condition obtained in Proposition \ref{prop:flux_condition} can be expressed with the help of the divergence-free in finite parts condition as a consequence of the following propositions. 
 
 \begin{proposition}\label{prop:div_implique_flux}
Let \(X:\R^d\rightarrow \R^d\) be a divergence-free in finite part vector field and let \(V\) be a smooth bounded open set. We assume that \(X\) is continuous in a neighbourhood  \(U\) of \(\p V\), then if \(\nu\) denotes the outer unit normal to \(\p V\),
\begin{equation*}
\int_{\p V} X\cdot \nu\dd \H^{d-1} =0.
\end{equation*}
\end{proposition}

\begin{proof}
Denote by $\textbf{1}_V$ the characteristic function of $V$. Let \((\zeta_n)_n\) be a sequence of functions in \(\C^\infty_c(\R^d)\) such that \(\zeta_n\rightarrow \textbf{1}_V\) in \(L^1(\R^d)\), \(0\leq 
\zeta_n\leq 1\) with \(\zeta_n=\textbf{1}_V\) outside \(U\) and \(0<\zeta_n<1\) in \(U\). We further assume about \((\zeta_n)_n\) that $|\nabla\zeta_n|\in\M_f^+(\R^d)$ and \(|\nabla \zeta_n|(\R^d)\rightarrow |D \textbf{1}_V|(\R^d)\) and finally that \(\nabla \zeta_n \rightharpoonup D \textbf{1}_{V}\) weakly-* in $(\C^0(\R^d;\R^d))^*$. For example we can take \(\zeta_n=\chi \rho_n \ast \textbf{1}_V\) for \(n\) large enough, where \(\chi\) is equal to \(1\) on \(V\setminus U\), \(\chi\) vanishes on \(\R^d \setminus (U\cup V)\) and \(0<\chi<1\) on \(U\). Denoting \(E=\R^d\setminus U\), we have \(\R^d \setminus \zeta_n^{-1}(\zeta_n(E))=\{x:\zeta_n (x)\neq 0,1\} \subset U\), and by the previous proposition we obtain 
\begin{equation*}
\int_{\R^d \setminus \zeta_n^{-1}(\zeta_n(E))} X\cdot \nabla \zeta_n\dd x=0.
\end{equation*}
On the other hand \(\nabla \zeta_n=0\) a.e.\ on \(\{x:\zeta_n(x)=0\}\) and \(\{x:\zeta_n(x)=1\}\), implying that \(\int_U X\cdot \nabla \zeta_n\dd x=0\). Then we pass to the limit  with \(n\to  +\infty\) and deduce that \(\int_U X\cdot D \textbf{1}_V\dd x= \int_{\p V} X\cdot \nu \dd \H^{d-1}=0.\)
\end{proof}

\begin{proposition}\label{prop:equiv_div_fp_zero_flux}
Suppose that \(X\in\C^0(\R^d\setminus \{x_1,\dots,x_N\};\R^d)\). Then \( \dive X=0\) in finite part, if and only if 
\begin{enumerate}
\item \(\dive X=0\) in \(\D'(\R^d\setminus \{x_1,\dots,x_N\})\)
\item \(\int_{\p B_\delta(x_i)} X\cdot \nu\dd \H^{d-1} =0\), for all \(i=1,\dots,N\) and all \(\delta< q_N=\inf\{|x_i-x_j|: i\neq j\}\).
\end{enumerate}
\end{proposition}

\begin{proof}
The first implication ``$\Rightarrow$'' is a consequence of Propositions \ref{prop:div_free_fp+L1_div_free} and \ref{prop:div_implique_flux}. We prove the converse implication ``$\Leftarrow$''. Take any \(\delta<q_N:= \min_{i\neq j} |x_i-x_j|\) and define 
\begin{equation*}
E_\delta :=\bigcup_{i=1}^N B_{\min(\frac{\delta}{N},q_N)} (x_i)
\end{equation*}
For \(\zeta \in \C^\infty_c(\R^d)\), by the coarea formula, we have
\begin{equation}\label{eq:coarea_2}
\int_{\R^d\setminus \zeta^{-1}(\zeta(E_\delta))} X\cdot \nabla \zeta\dd x=\int_{\{t:t \notin \zeta(E_\delta)\}} \int_{\{x:\zeta=t\}}X\cdot \nu \dd \H^{d-1} \dd t.
\end{equation}
Now by Sard's theorem for almost all \(t\in \R \setminus \zeta(E_\delta)\), the set \(\{x:\zeta=t\}\) is a smooth closed hyper-surface. We observe that this hyper-surface does not intersect any of the balls $\overline{B}_\delta(x_i)$, and thus
there is a finite number of balls (maybe zero) such that 
\( \overline{B}_\delta(x_i) \subset \{x:\zeta>t\}\), while the remaining balls are outside of its closure $\overline{\{x:\zeta>t\}}$. 
Without loss of generality we assume that the first \(N_p\) balls, with \(1\leq N_p\leq N\) are in  \(\{x:\zeta>t\}\). We apply the divergence theorem in \(\{x:\zeta>t\} \setminus \bigcup_{i=1}^{N_p} \overline{B}_\delta(x_i)\) and obtain 
\begin{align*}
\int_{\{x:\zeta=t\}}X\cdot \nu \dd \H^{d-1}&=\int_{\{x:\zeta>t\}\setminus \bigcup_{i=1}^{N_p} \overline{B}_\delta(x_i)} \dive X \dd x - \sum_{i=1}^{N_p}\int_{\p B_\delta(x_i)} X\cdot \nu \dd \H^{d-1}=0,
\end{align*}
where we have used  that \(\dive X=0\) in \(\D'(\R^d\setminus\{x_1,\dots,x_N\})\) and the zero flux condition in item (2) of the statement of the proposition.
Hence \eqref{eq:coarea_2} implies that \(\int_{\R^d\setminus \zeta^{-1}(\zeta(E_\delta))} X\cdot \nabla \zeta \dd x=0\). Since \(\mathcal{H}^1_\infty\left(\bigcup_{i=1}^N B_{\min(\delta/N,q_N)}(x_i)\right)\leq \delta\), the result follows.
\end{proof}

\begin{remark}
    \rm 
    The results in this section, in particular Proposition \ref{prop:equiv_div_fp_zero_flux} and Corollary \ref{cor:flux_condition_imp}, imply that \(\vx_N\in (\R^d)^N\) satisfies \eqref{eq:critical_points} if and only if \( [h_N,h_N]+R_N\) is divergence-free in finite part, where \(h_N\) is defined in \eqref{eq:3crit}.
\end{remark}

\section{Convergence of empirical measures and electric potentials}\label{sec:convergence}
In this section, we study the convergence, up to a subsequence, of empirical measures \(\mu_N\), and  of their associated electric potentials \(h_N\), cf.\ \eqref{eq:3crit}. In other words, we prove part A of Theorem \ref{th:main_1}. We begin with general results for the narrow and $W^{-1,p}$ convergence of finite signed Radon measures 
and for their associated electric potentials.

\begin{lemma}\label{lem:tightness}
    Let \((\mu_N)_N\) be a sequence of empirical measures \eqref{eq:3crit} satisfying boundedness condition \eqref{eq:Boundedness_Hamiltonian}. Then \( (\mu_N)_N\) is tight. 
\end{lemma}

\begin{proof}
    We denote \(A\wedge B:= \min (A,B)\). Consider the sequence of probability measures \((|\mu_N|)_N= \left(\frac{1}{M_N}\sum_{i=1}^N |d_i| \delta_{x_i}\right)_N\). 
     First note that, for any constant \(C_2>0\), there exists a compact set \(K \subset \R^d\) such that 
    \begin{equation}\label{eq:C_2}
        \min_{(K\times K)^c} \left[ g(x-y)+\frac{V(x)}{2}+\frac{V(y)}{2}\right] >C_2 .
    \end{equation}
    This follows from \cite[Lemma 2.10]{Serfaty_2015}. For example, when \(d=2\) we have 
    \begin{equation}\label{eq:g}
            g(x-y)=-\log |x-y| \geq -\log 2-\log \max(|x|,|y|),
    \end{equation}
    and thus \(\frac{1}{2}(V(x)+V(y))+g(x-y)\) is arbitrarily large if \(|x|\) and \(|y|\) are large enough. Therefore, if we rewrite
    \begin{align*}
        \frac{1}{2M_N^2}\sum_{i=1}^N \sum_{j\neq i} |d_i||d_j|\left[ g(x_i-x_j)+|F(x_i-x_j)|\right] +\frac{1}{M_N} \sum_{i=1}^N |d_i| V(x_i) \\
        = \iint_{\Delta^c} ( g(x-y)+|F(x-y)|) \dd |\mu_N|(x) \dd |\mu_N|(y) + \int_{\R^d} V(x) \dd |\mu_N|(x),
    \end{align*}
    where \(\Delta :=\{(x,x)\in \R^2: x\in \R\}\), it follows from the boundedness assumption \eqref{eq:Boundedness_Hamiltonian}, that there exists \(C_1>0\) such that, for every \(C>0\),
    \begin{align*}
        C_1 &\geq \iint_{\R^d\times \R^d} g(x-y)\wedge C + |F(x-y)| \dd |\mu_N|(x) \dd |\mu_N|(y) -\frac{C}{M_N} +\int_{\R^d} V(x) \dd |\mu_N|(x) \\
        & \geq \iint_{\R^d\times \R^d} \left( g(x-y)\wedge C+\frac{V(x)}{2}+\frac{V(y)}{2} \right) \dd |\mu_N|(x) \dd |\mu_N|(y) -\frac{C}{M_N}.
    \end{align*}
    Then growth properties of $V$, cf. \eqref{A2}, ensure existence of \(C_3>0\) such that \( g(x-y)+\frac{V(x)}{2}+\frac{V(y)}{2} \geq -C_3\) for any \(x,y\in \R^d\). Thus, using \eqref{eq:C_2}, we infer that for any \(C_2>0\), there exists a compact set \(K\subset \R^d\), such that
    \begin{align*}
         C_1\geq -C_3 -\frac{C}{M_N}+C_2 (|\mu_N| \otimes |\mu_N|)((K\times K)^c) \geq -C_3-\frac{C}{M_N}+C_2|\mu_N|(K^c)^2.
    \end{align*}
    Since \(C_2\) can be chosen arbitrarily large, \(|\mu_N|(K^c)\) can be made arbitrarily small and this means that \((|\mu_N|)_N\) is tight and, by definition, so is \((\mu_N)_N\) as a family of signed Radon measures. 
\end{proof}

\begin{proposition}\label{prop:compact}
    Let $(\mu_N)_{N\in\N}$ be a sequence in $\M_f(\R^d)$, which as a family of measures is bounded and tight. Then, up to a subsequence, $\mu_N\to \mu$ in $W^{-1,p}(\R^d)$ for some $\mu\in W^{-1,p}(\R^d)$ and for all $p\in[1, \frac{d}{d-1})$.
\end{proposition}

\begin{proof}
    By Prokhorov's theorem for signed measures, see e.g.\ \cite[Theorem 1.4.11]{Bogachev_2018}, $(\mu_N)_{N\in\N}$ is precompact in the narrow topology, and consequently there exists a signed Radon measure $\mu$, such that for all bounded-continuous functions $\phi:\R^d\to \R$, up to a subsequence, we have
    $$ \lim_{N\to +\infty}\int_{\R^d}\phi\dd \mu_N = \int_{\R^d}\phi \dd\mu.$$
    Exploiting tightness of $(\mu_N)_{N\in\N}$, we fix $\varepsilon>0$, and a compact set \(K\)
such that     
    $$ \sup_{N\in\N}|\mu_N|(\R^d\setminus K) <\varepsilon,\quad |\mu|(\R^d\setminus K)<\varepsilon. $$
 We can also take a bounded open set \(\Omega \subset \R^d\), such that \(K \subset \subset \Omega\). Denote \(\nu_N:= \mu_N \mres \Omega\) and \ $\nu:=\mu\mres \Omega$.  
Taking any $\phi\in \C_0^0(\Omega)$ and extending it
by \(0\) outside \(\Omega\) we find that
\begin{align}\label{eq:nuconv}
     \int_{\Omega} \phi \dd \nu_N &=\int_{\R^d}\phi \dd \mu_N
    \xrightarrow[]{N \to +\infty} \int_{\R^d} \phi \dd \mu = \int_{\Omega} \phi \dd \nu.
\end{align}

    \noindent
    Since \(\Omega\) is bounded we aim to use the Rellich-Kondrachov Theorem to prove that, up to a subsequence, \(\nu_N\) converges strongly to \(\nu\) in \(W^{-1,p}(\Omega)\) for any \(1\leq p<d/(d-1)\). Denote by $B$ the unit ball in $W_0^{1,p'}(\Omega)$. Our goal is to show that for $p'=\frac{p}{p-1}$,
    \begin{equation}\nonumber
    \|\nu_N - \nu\|_{W^{-1,p}(\Omega)}:= \sup_{\phi\in B}\left|\int_{\R^d}\phi\dd\nu_N - \int_{\R^d}\phi\dd\nu\right|\xrightarrow{N\to +\infty} 0.
    \end{equation}
    By the Rellich-Kondrachov Theorem, $B$ is compact in $\C^0_0(\Omega)$. Thus, given $\delta>0$, there exist functions $\{\phi_i\}_{i=1}^{n_\delta}$ in \(B\) such that for any $\phi\in B$ there exists $i\in\{1,...,n_\delta\}$ such that
    $$ \|\phi_i-\phi\|_{L^\infty(\Omega)}<\delta .$$
    Thus for any $\phi\in B$

    \begin{align*}
    \left|\int_{\Omega}\phi\dd\nu_N - \int_{\Omega}\phi\dd\nu\right| & \leq 2\|\phi\|_{L^\infty(\Omega)}\delta\sup_{N}|\nu_N|(\Omega) +  \left|\int_{\Omega}\phi_i\dd\nu_N - \int_{\Omega}\phi_i\dd\nu\right| 
    \end{align*}
    for some $i\in\{1,..., n_\delta\}$. In particular
    \begin{align*}
       \sup_{\phi \in B } \left|\int_{\Omega}\phi\dd\nu_N - \int_{\Omega}\phi\dd\nu\right| & \leq 2C \delta \sup_{N}|\nu_N|(\Omega) +  \max_{1\leq i \leq n_\delta}\left|\int_{\Omega}\phi_i\dd\nu_N - \int_{\Omega}\phi_i\dd\nu\right|,
    \end{align*}
    where here \(C\) is such that \(\|\phi\|_{L^\infty(\Omega)}\leq C \|\phi\|_{W^{1,p'}(\Omega)}\) thanks to the Sobolev embedding.
    Recalling that $\phi_i\in \C_0^0(\Omega)$, by \eqref{eq:nuconv}, we deduce that 
    \begin{align*}
        \limsup_{N \to +\infty} \sup_{\phi \in B } \left|\int_{\Omega}\phi\dd\nu_N - \int_{\Omega}\phi\dd\nu\right| \leq 2C \delta \sup_{N}|\nu_N|(\Omega).
    \end{align*}
    Since $\mu_N$ is bounded, we have $\sup_N|\nu_N|(\Omega)\leq \sup_N|\mu_N|(\Omega)< +\infty$, and since the last inequality holds for all \(\delta>0\) we conclude that \(\lim_{N \to +\infty} \| \nu_N-\nu\|_{W^{-1,p}(\Omega)}=0\).

    \medskip
Now let \(\chi\) be in \( \C^\infty_c(\Om)\) with \(\chi \equiv 1\) on \(K\) and \(0\leq \chi\leq 1\). For any \(\phi \in W^{1,p'}(\R^d) \) we  have \(\chi \phi \in W^{1,p'}_0(\Omega)\) and \(\| \chi \phi\|_{W^{1,p'}(\Om)}\leq C_{\chi} \|\phi\|_{W^{1,p'}(\R^d)}\). Since \(\phi =\chi \phi+(1-\chi) \phi\) we see that
\begin{align*}
    \left|\int_{\R^d}\phi\dd\mu_N - \int_{\R^d}\phi\dd\mu\right| &\leq \left| \int_{\Omega} \chi \phi \dd \mu_N -\int_{\Omega} \chi \phi \dd \mu \right| +2 \|\phi\|_{L^{\infty}(\R^d)} \varepsilon.
\end{align*}
    Thus, with \(C\) which denotes the constant in the Sobolev embedding \(W^{1,p'}(\R^d) \hookrightarrow \C^0_0(\R^d)\),
    \begin{align*}
          \left|\int_{\R^d}\phi\dd\mu_N - \int_{\R^d}\phi\dd\mu\right| \leq C_{\chi} \|\phi\|_{W^{1,p'}(\R^d)}\|\nu_N-\nu\|_{W^{-1,p}(\Omega)}+2C\varepsilon \|\phi\|_{W^{1,p'}(\R^d)}.
    \end{align*}
    Thus
    \begin{equation}\nonumber
        \|\mu_N-\mu\|_{W^{-1,p}(\R^d)}\leq C_{\chi} \|\nu_N-\nu\|_{W^{-1,p}(\Omega)}+2C\varepsilon
    \end{equation}
    and
    \begin{equation}\nonumber
        \limsup_{N \to +\infty}\|\mu_N-\mu\|_{W^{-1,p}(\R^d)} \leq 2C\varepsilon.
    \end{equation}
     Since the above is valid for any \(\varepsilon\), we obtain the conclusion.
    \end{proof}

\noindent
The following proposition has an important purpose. Together with Lemma \ref{lem:tightness} and Proposition \ref{prop:compact}, it ensures that up to a subsequence, the gradients of electric potentials $\nabla h_N:= \nabla g*\mu_N$ converge in $L^p_{\rm loc}(\R^d;\R^d)$ to $\nabla g*\mu$, where $\mu$ is the narrow limit of the empirical measures $\mu_N$.

\begin{proposition}\label{prop:hn}
Let $1\leq p<\frac{d}{d-1}$ and $R>0$. For every measure $\mu$ in \(\M_f(\R^d)\) the convolution $\nabla g\ast \mu$ is well defined as an $L^p_{\rm loc}(\R^d;\R^d)$ function and there exists $C_R>0$, depending only on $R$, $p$ and $d$, such that
\begin{equation}\label{e:hn0}
    \|\nabla g\ast \mu\|_{L^p(B_R(0))}\leq C_R\|\mu\|_{W^{-1,p}(\R^d)}.
\end{equation}
\end{proposition}

\begin{proof}
To prove that $\nabla g\ast\mu$ belongs to $L^1_{\rm loc}(\R^d;\R^d)$ it suffices to show that for any $R>0$ the convolution $|\nabla g|\ast |\mu|$ belongs to $L^1(B_R(0))$, where $|\mu|$ is the variation measure of $\mu$. By Fubini's theorem, we have
 \begin{equation*}
 \begin{split}
 \int_{B_R(0)} |\nabla g| \ast |\mu|(x) \dd x &=\int_{B_R(0)}\int_{\R^d} |\nabla g(x-y)|\dd |\mu|(y)\dd x\\
 &= \int_{\R^d}\int_{B_R(0)}|\nabla g(x-y)|\dd x\dd|\mu|(y)\\
 &= \int_{\R^d}\int_{B_R(y)}|\nabla g(z)|\dd z\dd|\mu|(y).
 \end{split}
 \end{equation*}
 Since $|\nabla g| \in L^1_{\text{loc}}(\R^d)$ and $|\nabla g(x)| \to 0$ as $|x| \to  +\infty$, it can be shown that the function $y \mapsto \int_{B_R(y)} |\nabla g(z)| \dd z$ is continuous and bounded. Since $|\mu|$ is finite, we find that $\int_{B_R(0)} |\nabla g| \ast |\mu|(x) \dd x \lesssim |\mu|(\R^d)$, and thus $\nabla g\ast\mu$ belongs to $L^1_{\rm loc}(\R^d;\R^d)$.

    To prove that actually $\nabla g\ast\mu\in L^p_{\rm loc}(\R^d;\R^d)$ it suffices to show \eqref{e:hn0}. Fix $R>0$ and any $\phi\in \C^\infty_c(\R^d)$ with $\|\phi\|_{L^{p'}(B_R(0))}\leq 1$. Fix \(1\leq i\leq d\) in the rest of the proof. Since $\p_i g\ast\mu\in L^1_{\rm loc}(\R^d)$, the function $\phi\p_i g\ast\mu$ belongs to $L^1(\R^d)$. Thus, given $\varepsilon>0$,   after possibly modifying $\phi$ outside of $B_R(0)$, we can assume that $\supp\phi\subset B_{2R}(0)$ and
    $$ \left|\int_{B^c_R(0)}\phi\p_i g\ast\mu\dd x \right| + \|\phi\|_{L^{p'}(B_R^c(0))}\leq \varepsilon. $$

    \noindent
    Then, by Fubini's theorem, we have
    \begin{equation}\label{eq:ixix}
             \begin{split}
                  \int_{B_R(0)}\phi \p_ig*\mu \dd x &=\int_{\R^d}\phi \p_i g*\mu \dd x + {\mathcal L}\qquad (|{\mathcal L}|\leq \varepsilon)\\
             &=\int_{\R^{d}}\int_{\R^d}\phi(x)\p_i g(x-y)\dd x\dd {\mu}(y)+{\mathcal L}\\
             &=\int_{\R^d}\int_{B_{2R}(0)}\p_i g(y-x) \phi(x)\dd x \dd\mu(y)+{\mathcal L},
             \end{split}
    \end{equation}
    and we turn our attention to the function
    $$ y\mapsto \int_{B_{2R}(0)}\p_ig(y-x) \phi(x)\dd x = \p_i g * \phi\quad (\mbox{since }\ \supp\phi\subset B_{2R}(0)).$$
    The goal is to bound it in $W^{1,p'}(\R^d)$ so that the finite signed measure $\mu$ can be bounded in $W^{-1,p}(\R^d)$ by duality.

    \noindent
    By Calderon-Zygmund inequality \cite[Theorem 9.9]{Gilbarg_Trudinger_2001} we have
    \begin{equation}\label{e:hn1} \|\nabla^2 (g\ast \phi)\|_{L^{p'}(\R^d)} \lesssim \|\Delta (g\ast \phi)\|_{L^{p'}(\R^d)} \leq \| \phi\|_{L^{p'}(B_R(0))} +\varepsilon,
    \end{equation}
    where $\nabla^2$ denotes the Hessian.

    \noindent
    On the other hand 
    $$ \int_{B_{2R}(0)}\nabla g(y-x) \phi(x)\dd x = \int_{B_{2R}(y)}\nabla g(z) \phi(z-y)\dd z , $$
    and thus, we shall consider
    \begin{equation*}
        \begin{split}
            \int_{\R^d}\left|\int_{B_{2R}(y)}\nabla g(z)\ \phi(z-y)\dd z\right|^{p'}\dd y &= \left(\int_{B_{4R}(0)}+\int_{ B_{4R}^c(0)}\right)\left|\int_{B_{2R}(y)}\nabla g(z)\ \phi(z-y)\dd z\right|^{p'}\dd y\\
    &=:I + II
        \end{split}
    \end{equation*}
    First, for $I$, using H\" older's inequality with $\frac{1}{p}+\frac{1}{p'}=1,$ we write
    $$ I\leq \int_{B_{4R}(0)}\|\phi\|_{L^{p'}(\R^d)}^{p'}\|\nabla g\|_{L^{p}(B_{2R}(y))}^{p'}\dd y. $$
    Since $\nabla g \in L^p_{\text{loc}}(\R^d;\R^d)$ and $\nabla g(x) \to 0$ as $|x| \to  +\infty$, it can be shown that the function \(y \mapsto \int_{B_{2R}(y)} |\nabla g(z)|^{p} \dd z\) is continuous and bounded, and thus
    \begin{equation}\label{e:hn2}
        I\lesssim \|\phi\|_{L^{p'}(\R^d)}^{p'} \leq (\|\phi\|_{L^{p'}(B_R(0))}+\varepsilon)^{p'}.
    \end{equation} 
 
\noindent
    To bound $II$, first note that for $y\notin B_{4R}(0)$ and $z\in B_{2R}(y)$, we have
    $$ |z|\geq \frac{|y|}{2}\quad \Leftrightarrow\quad |\nabla g(z)| \lesssim |z|^{1-d}\lesssim |y|^{1-d}  $$
    hence, using H\" older inequality again,
    \begin{equation}\label{e:hn3}
        II\lesssim \int_{B_{4R}^c(0)}|y|^{p'(1-d)} \|\phi\|_{L^{p'}(\R^d)}^{p'}\dd y \lesssim  (\|\phi\|_{L^{p'}(B_R(0))}+\varepsilon)^{p'}.
    \end{equation} 
    Here, the last inequality reflects the integrability of the function $y\mapsto |y|^{p'(1-d)}$ outside of $B_{4R}(0)$. Indeed, since $p<\frac{d}{d-1}$, we have $p'>d$, and thus $p'(1-d)< d-d^2\leq -d$ for $d\geq 2$. 

\noindent
    Combining inequalities \eqref{e:hn1}, \eqref{e:hn2} and \eqref{e:hn3} we infer that $$\|\p_i g\ast\phi\|_{W^{1,p'}(\R^d)}\lesssim \|\phi\|_{L^{p'}(B_R(0))} + \varepsilon,$$ which together with \eqref{eq:ixix} implies that
    $$ \left|\int_{B_R(0)} \phi\nabla g\ast \mu\dd x\right|\leq C_R(\|\phi\|_{L^{p'}(B_R(0))}+\varepsilon)\|\mu\|_{W^{-1,p}(\R^d)} +\varepsilon $$
    for all  $\phi\in \C^\infty(B_R(0))$ with $\|\phi\|_{L^{p'}(B_R(0))}\leq 1$ and for all $\varepsilon>0$. Therefore, taking supremum with respect to $\phi$, we obtain
    $$ \|\p_i g\ast\mu\|_{L^p(B_R(0))}\leq C_R\|\mu\|_{W^{-1,p}(\R^d)} $$
    and the proof is finished.
\end{proof}

\noindent
Finally, we prove the convergence of the reminders $R_N$, which together with results from Lemma \ref{lem:tightness} and Propositions \ref{prop:compact} and \ref{prop:hn} constitutes the proof of part A of Theorem \ref{th:main_1}. 

\begin{proposition}\label{prop:def_R}
    Let \(R_N\in L^p_{\loc}(\R^d;\R^{d\times d})\) be the regular reminder defined by \eqref{def:R_N}. Then there exists \(R\in L^p_{\loc}(\R^d; \R^{d\times d})\) such that, up to a subsequence,
    \begin{align}
        R_N &\rightarrow R\qquad \text{ in } L^p_{\loc}(\R^d;\R^{d\times d}),\label{eq:conv_R_N}\\
        \dive R&=2c_d ( (\nabla F)_\mu +\mu \nabla V),\nonumber
    \end{align}
    where \(\mu\in \mathcal{M}_f(\R^d)\) is the limit of a subsequence of the empirical measure \((\mu_N)_N\), cf. Proposition \ref{prop:compact}, and \( (\nabla F)_\mu\) is defined in \eqref{eq:false_conv}.
\end{proposition}

\begin{proof}
   
   Since \(\nabla V \in \C^0(\R^d;\R^d)\), narrow convergence of $\mu_N\to\mu$, implies the narrow convergence
        \( \mu_N \nabla V \to \mu \nabla V\). Similarly, since \(\nabla F \in \C^0(\R^d;\R^d)\) and $\mu_N\otimes\mu_N\to \mu\otimes\mu$ narrowly (narrow convergence of the product measure), we  see that, for any \(\varphi \in \C^\infty_c(\R^d)\),
        \begin{multline}\label{eq:conv_F_mu}
            \iint_{\R^d \times \R^d}  \nabla F (x-y)\cdot (\varphi(x)-\varphi(y)) (\dd \mu_N\otimes \dd \mu_N)(x,y) \\
            \rightarrow\iint_{\R^d \times \R^d}  \nabla F(x-y)\cdot (\varphi(x)-\varphi(y)) (\dd \mu\otimes \dd \mu)(x,y).
        \end{multline}
        By using the symmetry of \(F\) and Fubini's theorem, we also have 
        \begin{align}
            \langle \mu_N \nabla F \ast \mu_N, \varphi \rangle &= \int_{\R^d} \left( \nabla F(x-y) \dd \mu_N(y) \right) \cdot \varphi(x)\dd \mu_N(x)\nonumber \\
            & =\frac12 \iint_{\R^d \times \R^d} \nabla F(x-y) \cdot (\varphi(x)-\varphi(y))  (\dd \mu_N\otimes \dd \mu_N)(x,y). \label{eq:conv_F_mu_2}
        \end{align}
        Then \eqref{eq:conv_F_mu} and \eqref{eq:conv_F_mu_2} show that \(\mu_N \nabla F \ast \mu_N\) converges weakly-* towards \( (\nabla F)_\mu\),  defined in \eqref{eq:false_conv}.
        Now, using again Proposition \ref{prop:compact} and Proposition \ref{prop:hn} we infer that \((R_N)_{ik}^T \rightarrow \frac{-2\p_ig}{c_d}\ast [( (\p_k F)_\mu+\mu \p_k V)]=:R_{ik}^T\) in \(L^p_{\text{loc}}(\R^d)\) for \(1\leq p<d/(d-1)\). Here \((\p_k F)_\mu\) is the \(k\)-th component of \( (\nabla F)_\mu\).  Hence \(R_N \rightarrow R\) in \(L^p_{\text{loc}}(\R^d;\R^{d\times d})\) for \(1\leq p<d/(d-1)\) and, from the definition of \(R_N\) and the fact that \(g\) satisfies \(-\Delta g =c_d \delta_0\) we find that \(\dive R= 2 ((\nabla F)_\mu +\mu\nabla V)\).

\end{proof}

\section{ Criticality condition in the weak velocity formulation}\label{sec:velocity_form}

This section is devoted to the proof of part B of Theorem \ref{th:main_1}. The key problem is the construction of the family $\{E_\delta\}_{\delta>0}$ appearing in the divergence-free in finite part condition, cf. Definition \ref{def:div_free_fp}. We shall perform such a construction along the proofs of the following two propositions.

\begin{proposition}\label{prop:conv_capacity}
    Let $(\alpha_N)_N$ be a sequence in $\M_f(\R^d)$ converging to $0$ in $W^{-1,p}(\R^d)$ with $1\leq p<\frac{d}{d-1}$. Then there exists a (non-relabelled) subsequence  and a collection of sets $\{G_\delta\}_{\delta>0}$, such that for all compact subsets $K\subset\R^d$, we have
    $$ \capp_p(G_\delta\cap K)\xrightarrow{\delta\to 0} 0 $$
    and
    $$ \forall \delta>0\quad \int_{K\setminus G_\delta}| \nabla g \ast \alpha_N|^2\dd x\xrightarrow{k\to +\infty} 0.$$
\end{proposition}

\begin{proof}
    First note that it suffices to prove the result by taking all balls $B_\ell(0)$ with $\ell=1,2,...$ instead of arbitrary compact sets $K$. Indeed, since the capacity is an outer measure, if $K\subset B_\ell(0)$ then $\capp_p(G_\delta\cap K)\leq \capp_p(G_\delta\cap B_\ell(0))$ and 
    $$ \forall\delta>0 \quad \int_{K\setminus G_\delta}|\nabla g \ast \alpha_N|^2\dd x \leq \int_{B_\ell(0)\setminus G_\delta}| \nabla g \ast \alpha_N|^2\dd x. $$
    Let us construct the desired family $\{G_\delta\}_{\delta>0}$. To this end fix \(1\leq p<\frac{d}{d-1}\), $\ell\in \N$ and  $\delta>0$.  Let \(f_{N,\ell}\) be such that \(-\Delta f_{N,\ell}= \alpha_N\) in \(B_{\ell}(0)\) with \(f_{N,\ell}\in W^{1,p}_0(B_\ell(0))\)  and 
    $$\|f_{N,\ell}\|_{W^{1,p}_0(B_\ell(0))}\leq C_\ell \|\alpha_N\|_{W^{-1,p}(\R^d)}.$$
 Such an \(f_{N,\ell}\) exists, and is unique, we refer e.g.\ to \cite[Proposition 3.2 and Proposition 5.1]{Ponce_2016}. We extend \(f_{N,\ell}\) to a Sobolev function defined in all \(\R^d\), still denoted \(f_{N,\ell}\), and such that \(\|f_{N,\ell}\|_{W^{1,p}(\R^d)}\leq C_\ell \|\alpha_N\|_{W^{-1,p}(\R^d)} \).

    By \cite[Theorem 4.18 and Theorem 4.19]{Evans_Gariepy_2015}, each of the functions $f_{N,\ell}$ has a representative $f_{N,\ell}^*$ (equal up to a subset of $\capp_p$ zero to $f_{N,\ell}$) for which 
    \begin{equation}\label{e:cap1}
           \capp_p(\{x\in\R^d:\ f_{N,\ell}^*(x)>\eta\})\lesssim \frac{\|\nabla f_{N,\ell}\|_{L^p(B_\ell(0))}^p}{\eta^p}. 
    \end{equation}

  \noindent
    Since sets of capacity zero are null with respect to the Lebesgue measure (cf. \cite[Theorem 4.15]{Evans_Gariepy_2015}), throughout the proof we shall assume that $f_{N,\ell}=f^*_{N,\ell}$, and henceforth we forfeit the star superscript. Let
    $$ F_{N,\ell}:= \{x\in\R^d:\ |f_{N,\ell}|\geq \delta_N\},\quad \delta_N:=\left(\frac{\|\alpha_N\|_{W^{-1,p}}}{1+\|\alpha_N\|_{TV}}\right)^\frac{1}{2}. $$
    Then $\delta_N\to 0$ as $N\to +\infty$ and thus, up to a subsequence $F_{N,\ell}$ is an increasing collection with \(N\) (recall that $\ell$ is fixed). Moreover, by \eqref{e:cap1} and by Proposition \ref{prop:hn}, we have
    $$ \capp_p(F_{N,\ell})\leq C_\ell\frac{\|\nabla f_{N,\ell}\|_{L^p(B_\ell(0))}^p}{\delta_N^p}\leq \frac{C_\ell\|\alpha_N\|_{W^{-1,p}}^p}{\delta_N^p}\leq C_\ell\|\alpha_N\|_{W^{-1,p}}^\frac{p}{2}(1+\|\alpha_N\|_{TV})^\frac{p}{2}\xrightarrow{N\to +\infty} 0. $$
    Thus, for any $\ell$ there exists $N_\ell$ such that for all $N\geq N_\ell$, we have
    $$ \capp_p(F_{N,\ell})\leq \frac{1}{2^\ell},$$
    and consequently
    $$ \sum_{\ell=1}^{ +\infty}\capp_p(F_{N_\ell,\ell})< +\infty. $$
    The above ensures that
    $$ G_{\delta}:= \bigcup_{\ell>\frac{1}{\delta}}F_{N_\ell,\ell} $$
    satisfies
    $$ \capp_p(G_{\delta})\leq \sum_{\ell>\frac{1}{\delta}}\capp_p(F_{N_\ell,\ell})\xrightarrow{\delta\to 0} 0. $$
    It remains to show that for any fixed $\delta>0$ and $\ell_0\in\N$ , we have
    $$ \int_{B_{\ell_0}(0)\setminus G_\delta}|\nabla f_{N_\ell,\ell}|^2 \dd x=\int_{B_{\ell_0}(0)\setminus G_\delta}|\nabla g \ast \alpha_{N_\ell}|^2 \dd x\xrightarrow{\ell\to +\infty} 0. $$
    Let $\bar{f}_{N_\ell,\ell} = {\rm sgn} \ f_{N_\ell,\ell}\min\{|f_{N_\ell,\ell}|,\delta_N\}\in W^{1,p}\cap L^\infty(B_{\ell_0}(0))$ and consider for $\ell\geq \ell_0$
    \begin{align*} \int_{B_{\ell_0}(0)\setminus G_\delta}|\nabla f_{N_\ell,\ell}|^2\dd x &\leq \int_{B_{\ell_0}(0)\setminus F_{N_\ell,\ell}}|\nabla f_{N_\ell,\ell}|^2\dd x = \int_{B_{\ell_0}(0)}\nabla f_{N_\ell,\ell}\cdot\nabla \bar{f}_{N_\ell,\ell}\dd x \\
    &= -\int_{B_\ell(0)} \alpha_{N_\ell} \bar{f}_{N_\ell,\ell}\dd x \leq \|\alpha_{N_\ell}\|_{TV}\delta_{N_\ell}\to 0. 
    \end{align*}
    To justify the last equalities we need to approximate \(f_{N_\ell,\ell}\) by smooth functions with compact supports in \(B_\ell(0)\) called \( (f_\e)_\e\). As usual, such an approximation is obtained by convolution and truncation (after extending the function by zero outside \(B_\ell(0)\)), and we have \(\|\Delta f_\e\|\leq \|\alpha_{N_\ell}\|\). Moreover \(\bar{f}_\e:={\rm sgn} \ f_\e \min\{|f_\e|,\delta_N\} \) is Lipschitz with compact support and \(\nabla f_\e \rightarrow \nabla f_{N_\ell,\ell}\) a.e.\ and \(\nabla \bar{f}_\e \rightarrow \nabla\bar{f}_{N_\ell,\ell}\) a.e. Hence, by Fatou's Lemma,
    \begin{align*}
        \int_{B_\ell(0)} \nabla f_{N_\ell,\ell}\cdot \nabla \bar{f}_{N_\ell,\ell}\dd x&\leq \liminf_{\e \to 0} \int_{B_\ell(0)} \nabla f_\e \cdot \nabla \bar{f}_\e \dd x =-\int_{B_\ell(0)} \Delta f_\e \bar{f}_\e\dd x \leq \|\alpha_{N_\ell,\ell}\| \delta_N.
    \end{align*}

\end{proof}

Before we proceed with the construction of the collection $(E_\delta)_\delta$ and with the culmination of the proof of part B of Theorem \ref{th:main_1}, let us collect all the information provided by Propositions \ref{prop:compact}, \ref{prop:hn} and \ref{prop:conv_capacity}.

\begin{corollary}\label{cor:summary}
    Let \(\mu_N\) be the sequence of empirical measures associated with critical points of $\H_N$, and \(h_N\) be their electric potentials, cf. \eqref{eq:3crit}.  Then there exists \(\mu\in \mathcal{M}_f(\R^d)\) such that, up to a subsequence, for any \(1\leq p< d/(d-1)\),
    \begin{equation}\label{eq:conv_mu_mutilde}
        \mu_N \to \mu, \ \text{ narrowly in } \mathcal{M}_f(\R^d), \quad \nabla h_N \rightarrow \nabla g \ast \mu, \text{ in } W^{1,p}_{{\rm loc}}(\R^d;\R^d).
    \end{equation}
    Furthermore, there exists a family of set \(\{G_\delta\}_{\delta>0}\), such that for all compact subsets $K\subset\R^d$, we have
    \begin{equation}\label{eq:L1}
    \capp_p(G_\delta\cap K)\xrightarrow{\delta\to 0} 0, \text{ and }
      \int_{K\setminus G_\delta}| \nabla h_N- \nabla g \ast \mu|^2\dd x\xrightarrow{N\to +\infty} 0.
    \end{equation}
\end{corollary}

\begin{proposition}\label{prop:div_free_cap}
    Let \(\vx_N = (x_1,\dots,x_N)\in (\R^d)^N\) satisfy \eqref{eq:critical_points}. Let \(\mu_N,h_N\) be the sequence of empirical measures and their associated electric potentials defined in \eqref{eq:3crit}. Moreover, let $\mu$ be the limit ensured by Proposition \ref{prop:compact}. Then there exists a family of sets \(\{E_\delta\}_{\delta>0}\), such that for any compact set \(K \subset \R^d\), we have \(\lim_{\delta \to 0} \capp_p(K\cap E_\delta)=0\) for any \(1\leq p<d/(d-1)\) and for all \(\delta>0\) and all \(i=1,\dots,d\) 
    \begin{equation}\label{eq:weak_div_free}
        \int_{\R^d \setminus \zeta^{-1}(\zeta (E_\delta))} (\lh\nabla g*\mu, \nabla g*\mu\rh^i+R^i)\cdot \nabla \zeta \dd x=0 \quad \text{for all } \zeta\in \C^\infty_c(\R^d)
    \end{equation}
  where \(\lh\nabla g*\mu, \nabla g*\mu\rh^i+R^i\)  denotes the ith row of \(\lh\nabla g*\mu, \nabla g*\mu\rh+R\), with the reminder \(R\)  defined in Proposition \ref{prop:def_R}.
\end{proposition}

\begin{proof}
Let \(\zeta \in \C^\infty_c(\R^d)\) and let \(K:= \supp \zeta\). We take \( \{G_\delta\}_{\delta>0}\) the family of sets given by Corollary \ref{cor:summary}. Let $\{N_\ell\}_{\ell=1}^{ +\infty}$ be the subsequence introduced in the construction of $G_\delta$ in the proof of Proposition \ref{prop:conv_capacity}.

For each $\ell\in\N$, consider the support $\{x^{N_\ell}_1,..., x^{N_\ell}_{N_\ell}\}$ of the empirical measure $\mu^{N_\ell}$ and arrange all of these points in a countable set as
$$ x_1 = x^{N_1}_1,\ x_2 = x^{N_1}_2,...,\ x_{N_1} = x^{N_1}_{N_1},\ x_{N_1+1} = x^{N_2}_1,...$$
i.e. we list all of the supports of $\mu^{N_\ell}$ one after another starting with $\ell=1$. We define

\begin{equation}\label{gdelta}
    E_\delta:= G_\delta \cup  \bigcup_{i=1}^{ +\infty} B(x_i, (\delta 2^{-i})^\frac{1}{d-p}).
\end{equation} 
Since $\capp_p$ is an outer measure satisfying inequality $\capp_p(B(x,r))\leq r^{d-p}\capp_p(B(0,1))$ (cf. \cite[Theorem 4.15]{Evans_Gariepy_2015}), we have
$$ \capp_p(K\cap E_\delta)\leq \capp_p(K\cap G_\delta) + \delta\capp_p(B(0,1))\xrightarrow{\delta\to 0}0 .$$

Note that \( \lh \nabla h_N, \nabla h_N\rh+R_N\) and \(\lh\nabla g*\mu, \nabla g*\mu\rh+R\) are in \(L^1(K \setminus E_\delta) \) since \(\nabla h_N\) is in \(L^2(K\setminus E_\delta)\) and thanks to \eqref{eq:L1}.  Thus, for all \(i=1,\dots,d\),
\begin{align*}
    \Bigg|\int_{\R^d \setminus \zeta^{-1}(\zeta(E_\delta))} & (\lh \nabla h_N, \nabla h_N\rh^i+R_N^i- (\lh\nabla g*\mu, \nabla g*\mu\rh^i+R^i)) \cdot \nabla \zeta \dd x\Bigg| \\
    &\leq \int_{\{x \in K: \zeta (x) \notin \zeta(E_\delta)\}} \left|(\lh \nabla h_N, \nabla h_N\rh^i+R_N^i-(\lh\nabla g*\mu, \nabla g*\mu\rh^i + R^i))\right||\nabla \zeta | \dd x\\
    & \leq \| \nabla \zeta \|_{L^\infty} \int_{K \setminus E_\delta} \left|(\lh \nabla h_N, \nabla h_N\rh^i+R_N^i-(\lh\nabla g*\mu, \nabla g*\mu\rh^i + R^i))\right|\dd x
    \xrightarrow{N \to 0} 0,
\end{align*}
thanks to \eqref{eq:L1} and \eqref{eq:conv_R_N}. 

Then, denoting \(\gamma_t:= \{x\in K: \zeta(x)=t\}\), coarea formula yields
\begin{align}
    \int_{\R^d\setminus \zeta^{-1}(\zeta (E_\delta))} (\lh \nabla h_N, \nabla h_N\rh^i+R_N^i)\cdot \nabla \zeta\dd x &= \int_{\{t:\gamma_t\cap E_\delta =\emptyset\}} \int_{\gamma_t} (\lh \nabla h_N, \nabla h_N\rh^i+R_N^i)\cdot \nu\dd\H^{d-1}. \label{eq:integral_t}
\end{align}

\noindent
However, for any regular value $t$ of \(\zeta\) satisfying \(\gamma_t\cap E_\delta =\emptyset\) we can apply the divergence theorem, or directly Corollary \ref{cor:flux_condition_imp}, to obtain 
\begin{equation*}
    \int_{\gamma_t} (\lh \nabla h_N, \nabla h_N\rh^i+R_N^i)\cdot \nu \dd\H^{d-1}=0.
\end{equation*}
Thus, in light of \eqref{eq:integral_t}, we have 
\begin{equation}\label{eq:zero_on small_capa}
     \int_{\R^d\setminus \zeta^{-1}(\zeta (E_\delta))} (\lh \nabla h_N, \nabla h_N\rh^i+R_N^i)\cdot \nabla \zeta\dd x=0. 
\end{equation}

Now we write
\begin{align*}
    \int_{\R^d \setminus \zeta^{-1}(\zeta (E_\delta))} &(\lh\nabla g*\mu, \nabla g*\mu\rh^i+R^i)\cdot \nabla \zeta\dd x \\
   & = \int_{\R^d \setminus \zeta^{-1}(\zeta (E_\delta))}(\lh\nabla g*\mu, \nabla g*\mu\rh^i+R^i-(\lh \nabla h_N, \nabla h_N\rh^i+R_N^i))\cdot \nabla \zeta\dd x\\
&\qquad +\int_{\R^d \setminus \zeta^{-1}(\zeta (E_\delta))}(\lh \nabla h_N, \nabla h_N\rh^i+R_N^i)\cdot \nabla \zeta\dd x.
\end{align*}
Since we have proved that the first term on the right-hand side above converges to zero as \(N \rightarrow  +\infty\) and that the second term is equal to zero,  
\eqref{eq:weak_div_free} follows.
\end{proof}

Observe that, by Corollary \ref{cor:summary} and Proposition \ref{prop:div_free_cap}, the vector field $\lh\nabla g*\mu, \nabla g*\mu\rh^i + R^i$, together with family $\{E_\delta\}_\delta$ satisfies the last two items of Definition \ref{def:div_free_fp}. Instead of the first item, all we can say is that
$$ \capp_p(K\cap E_\delta)\xrightarrow{\delta\to 0} 0. $$
To make sure that actually $\H^1_\infty(K\cap E_\delta)\to 0$ (which then implies that $\lh\nabla g*\mu, \nabla g*\mu\rh^i+ R^i$ is divergence-free in finite part), we use the following lemma.

\begin{lemma}\label{lem:importance_dim_2}
    There exists \(C>0\) such that for any bounded set \(A\subset \R^2\) we have 
    \begin{equation*}
        \mathcal{H}^1_\infty(A)\leq C \capp_1 (A). 
    \end{equation*}
    In particular, see \eqref{eq:Lipschitz_rule}, for any Lipschitz function \(\zeta:\R^d\to\R\) there holds \(\mathcal{L}^1( \zeta (A))\leq C \|\zeta\|_{\text{Lip}} \capp_1(A).\)
\end{lemma}

The above lemma is the only part of our result, where restriction to dimension 2 is crucial. Indeed, in dimension \(d\geq 3\) we cannot control the 1-Hausdorff content of a bounded set \(A\subset \R^d\) by its \(p\)-capacity for \( 1\leq p<d/(d-1)\). Lemma \ref{lem:importance_dim_2} is proved in \cite[Lemma 13.1]{Sandier_Serfaty_2007} and is a reformulation of \cite[Theorem 5.12]{Evans_Gariepy_2015}. 

\subsection{Proof of Theorem \ref{th:main_1} -- a summary}
Throughout Section \ref{sec:convergence} and Section \ref{sec:velocity_form} we proved Theorem \ref{th:main_1}. Let us perform a checklist of the proof to ensure that everything is in order. 
\begin{proof}[Proof of Theorem \ref{th:main_1}]
    We begin with the convergence in part A, which we proved in Section \ref{sec:convergence}. The fact that \(\mu_N \to \mu\) narrowly follows from Proposition \ref{prop:compact}. The fact that \(\nabla h_N \rightarrow \nabla g \ast \mu\) in \(L^p_{\text{loc}}(\R^d;\R^d)\) is a consequence of Proposition \ref{prop:compact} and Proposition \ref{prop:hn}. The convergence of \(R_N\) follows from Proposition \ref{prop:def_R}. Criticality condition in part B was dealt with in Section \ref{sec:velocity_form}, particularly Proposition \ref{prop:div_free_cap} and Lemma \ref{lem:importance_dim_2} ensure the divergence-free in finite part condition for the tensor $\lh\nabla g*\mu, \nabla g*\mu\rh^i + R^i$. The family $E_\delta$ was constructed in Propositions \ref{prop:conv_capacity} and \ref{prop:div_free_cap}.
\end{proof}

\section{Criticality condition in the velocity formulation}\label{sec:aug}

In this short section, we prove Theorem \ref{th:main1.5}, which is based on an additional assumption that the limit $\mu$, established in part A, belongs to $H^{-1}_{\rm loc}(\R^2)$ and leads to an augmented criticality condition of Definition \ref{def:critcond}.

\begin{proof}[Proof of Theorem \ref{th:main1.5}]
    First, observe that \eqref{eq:1.16} follows directly from Theorem \ref{th:main_1} and Proposition \ref{prop:div_free_fp+L1_div_free} in the preliminaries. This is because if $\mu\in H^{-1}_{\rm loc}(\R^2)$, then $\nabla g*\mu\in L^2_{\rm loc}(\R^2;\R^2)$ and thus the troublesome singular part of the considered vector field, $\lh\nabla g*\mu, \nabla g*\mu\rh^i$, belongs to $L^1_{\rm loc}(\R^2;\R^2)$.

    \medskip
    \noindent
    It remains to prove \eqref{eq:possible_reg}. We set \(Y:=\nabla g \ast \mu + \nabla F \ast \mu +\nabla V\), then
        \begin{align*}
            \lh Y,Y\rh &=\lh\nabla g \ast \mu, \nabla g \ast \mu\rh +\lh\nabla F \ast \mu, \nabla F \ast \mu\rh +\lh\nabla V, \nabla V\rh \\
            &+2 \left( \lh\nabla g \ast \mu, \nabla F \ast \mu\rh+\lh\nabla g \ast \mu, \nabla V\rh+ \lh\nabla F \ast \mu,V\rh \right).
        \end{align*}
        From \eqref{eq:1.16} we know that \( \dive \lh\nabla g \ast \mu, \nabla g \ast \mu\rh=-2\mu (\nabla F \ast \mu +\nabla V) \) in the sense of distributions. Moreover, by a minor modification of Lemma \ref{lem:chain_rule_stress_tensor}, the following identities hold true:
        \begin{align*}
        \dive \lh\nabla F \ast \mu,\nabla F\ast \mu \rh&=2\dive ( \nabla F \ast \mu )\nabla F \ast \mu, \\
        \dive \lh\nabla V, \nabla V\rh &= 2\Delta V \nabla V, \\
        \dive \lh\nabla g \ast \mu, \nabla F \ast \mu\rh&=\mu \nabla F \ast \mu +\dive ( \nabla F \ast \mu ) \nabla g \ast \mu, \\
        \dive \lh\nabla g \ast \mu, \nabla V\rh&=\mu \nabla V + \Delta V \nabla g \ast \mu, \\
        \dive \lh\nabla F \ast \mu,V\rh&= \dive (\nabla  F \ast \mu ) \nabla V + \Delta V \nabla F \ast \mu.
        \end{align*}
        Hence we find that 
        \begin{align*}
            \dive \lh Y,Y\rh& = -2\mu (\nabla F \ast \mu +\nabla V)+2 \dive ( \nabla F \ast \mu) \nabla F \ast \mu+2\Delta V \nabla V  \\
            & \quad +2(\mu \nabla F \ast \mu +\dive (\nabla  F \ast \mu) \nabla g \ast \mu)+2(\mu \nabla V + \Delta V \nabla g \ast \mu)  \\
            & \quad \quad +2(\dive (\nabla  F \ast \mu ) \nabla V + \Delta V \nabla F \ast \mu) \\
            &=2 (\dive (\nabla  F \ast \mu)+\Delta V) (\nabla g \ast \mu+\nabla F \ast \mu+\nabla V).
        \end{align*}

        \smallskip
        \noindent
        We now prove that \( \nabla g \ast \mu \in L^\infty_{\text{loc}}(\R^2;\R^d)\) if $\Delta F$ and $\Delta V$ belong to \(L^\infty_{\text{loc}}(\R^2)\). To do that, we take advantage of the complex structure of \(\R^2 \simeq \mathbb{C}\). We write \( Y=(Y_1,Y_2)\in \R^2\) and we identify it with \( \bar{Y}=Y_1-iY_2\in \mathbb{C}\). We recall that \(\p_{\bar{z}}=\frac{1}{2}(\p_x+i\p_{y})\) and hence
        \begin{equation}\label{partzet}
        \begin{split}
            2\p_{\bar{z}}[\bar{Y}^2]&=\left[ \p_x(Y_1^2-Y_2^2)+\p_y(2Y_1Y_2) \right] +i\bigl[ \p_y (Y_1^2-Y_2^2)-\p_1 (2Y_1Y_2)\bigr] \\
            &= (\dive [Y,Y])_1-i( \dive [Y,Y])_2 \\
            &= 2 (\Delta F \ast \mu + \Delta V) (Y_1-iY_2).
        \end{split}
        \end{equation}
        Hence, from the assumption that \(\mu \in H^{-1}_{\text{loc}}(\R^2)\), we see that \(\nabla g \ast \mu \in L^2_{\text{loc}}(\R^2;\R^2)\), and thus \(\p_{\bar{z}} [\bar{Y}^2]\in L^2_{\text{loc}}(\R^2;{\mathbb C})\). Since the operator \(\p_{\bar{z}}\) is elliptic (the function \(z\mapsto\frac{1}{z}\) is, up to a multiplicative constant, a fundamental solution for this operator) we conclude that \( \bar{Y}^2 \in W^{1,2}_{\text{loc}}(\R^2;{\mathbb C})\). Thanks to Sobolev's embedding we infer that \(\bar{Y}^2\) belongs to \(L^p_{\text{loc}}(\R^2;{\mathbb C})\) for any \(1\leq p < +\infty\), and so does \(\bar{Y}\). We bootstrap this information using \eqref{partzet}, because now we know that \(\p_{\bar{z}}[\bar{Y}^2] \in L^p_{\text{loc}}(\R^2;{\mathbb C})\) for any \(1\leq p< +\infty\). Thus \(\bar{Y}^2\in W^{1,p}_{\text{loc}}(\R^2;{\mathbb C})\) for any \(1\leq p< +\infty\), and by Sobolev's embedding, \(|Y|^2 \in \C^{0,\alpha}(\R^2)\) for any \(0<\alpha<1.\) By taking the square root we conclude that \(|Y| \in \C^{0,\alpha}(\R^2)\) for any \(0\leq \alpha < 1/2\) and it is locally bounded. Since \(\nabla g \ast \mu=Y-\nabla F \ast \mu -\nabla V\), and since \(\nabla F\ast \mu\) and \(\nabla V \) are continuous, it follows that \(\nabla g \ast \mu \in L^\infty_{\text{loc}}(\R^2;\R^2)\). Now assuming \(\mu \in L^1_{\text{loc}}(\R^2)\), and keeping in mind that \(\nabla g \ast \mu \in L^\infty_{\text{loc}}(\R^2;\R^2)\), we apply Lemma \ref{lem:chain_rule_stress_tensor} to obtain \(\dive [\nabla g \ast \mu, \nabla g \ast \mu]=2\mu \nabla g \ast \mu\). Hence from item (1) the following stronger form  is obtained:
        \begin{equation}\label{strongerform}
            \mu ( \nabla g \ast \mu+\nabla F \ast \mu +\nabla V)=0.
        \end{equation}

       \noindent
        Finally, to prove \eqref{eq:description} we observe that under the assumption \(\mu \in L^1(\R^2)\), equation \eqref{strongerform}  implies that \(\mu=0\) a.e.\ on \(\{x:|\nabla g \ast \mu+\nabla F \ast \mu+\nabla V|>0\}\). But we can use  \cite[Theorem 1]{Ambrosio_Ponce_Rodiac_2020} to infer that \(\Delta g \ast \mu +\Delta (F \ast \mu)+ \Delta V =0\) a.e.\ on \(\{x:|\nabla g \ast \mu+\nabla F \ast \mu+\nabla V|=0\}\). Since \(-\Delta g \ast \mu=c_d \mu\) we find that 
        $$\mu=\frac{(\Delta (F \ast \mu)+ \Delta V)\textbf{1}_{\{|\nabla g \ast \mu+\nabla F \ast \mu+\nabla V|=0\}}}{c_d}.$$
\end{proof}

\section{Criticality condition in the vorticity formulation}\label{sec:vorticity_form}

We dedicate the following section to the proof of Theorem \ref{th:main_2}; that is we pass to the limit in the criticality conditions satisfied by the empirical measures \(\mu_N\) in the vorticity form. As before we assume that \(F\in \C^1\).
We begin with the following lemma, whose proof can be found for instance in \cite[Lemma 5.2]{Schochet_1995}.

\begin{lemma}\label{lem:domination_measure}
    Let \( (\alpha_N)_N\subset \mathcal{M}_f(\R^d)\) and \( (\beta_N)_N \subset \mathcal{M}_f^+(\R^d)\) be such that there exists \(C>0\) with \(|\alpha_N|\leq C\beta_N\) for all \(N\in \mathbb{N}\). We assume that \(\alpha_N \to \alpha \)  and \(\beta_N \to \beta\) narrowly. Then there exists a \(\beta\)-measurable function \(f:\R^d \rightarrow \R\)  such that \(\alpha=f\beta\) and such that \(|f(x)|\leq C\) for \(\beta\)-a.e. \(x\).
\end{lemma}

 The following is inspired by \cite[Lemma 5.3]{Schochet_1995}.

\begin{lemma}\label{lem:factorization}
\noindent
For all \(\varphi \in \C^\infty_c(\R^2;\R^2)\) we have
\begin{equation}
    \nabla g (x-y)\cdot  \left(\varphi (x)-\varphi(y) \right)=Q_1(x-y) G_{1,\varphi}(x,y)+Q_2(x-y)G_{2,\varphi}(x,y)+Q_3(x-y)G_{3,\varphi}
\end{equation}
    where \(Q_i\in \C^0(\R^2\setminus \{0\})\), \(i=1,2,3\) are bounded and independent of \(\varphi\) and \(G_{i,\varphi}(x,y) \in \C^0_0(\R^2\times \R^2)\).
\end{lemma}

\begin{proof}
    We use the fundamental theorem of calculus to write
    \begin{align*}
        \varphi(x)-\varphi(y) &=\int_0^1 \frac{\dd}{\dd s}\varphi(sx +(1-s)y) \dd s =\left[\int_0^1 D\varphi (s x +(1-s)y) \dd s\right] (x-y),
    \end{align*}
    where we write $[A]b$ above to indicate that $[A]$ is a matrix acting on the vector $b$.
Hence
    \begin{align*}
        \nabla g (x-y)\cdot\left(\varphi (x)-\varphi(y) \right)  &= (\nabla g (x-y))^T\left[\int_0^1 D \varphi (sx +(1-s)y) \dd s\right] (x-y)  \\
         &=\sum_{k=1}^2 \sum_{i=1}^2  \int_0^1 \p_k\varphi^i (sx+(1-s)y) \dd s (x-y)_k \p_i g (x-y) \\
        &= \frac{-(x_1-y_1)^2}{|x-y|^2} \int_0^1 \p_1 \varphi^1  (sx+(1-s)y) \dd s \\
        & \quad -\frac{(x_2-y_2)^2}{|x-y|^2}\int_0^1 \p_2 \varphi^2  (sx+(1-s)y) \dd s \\
        & \quad \quad  -\frac{(x_1-y_1)(x_2-y_2)}{|x-y|^2} \Bigl( \int_0^1 \p_1 \varphi^2  (sx+(1-s)y) \\
        & \quad \quad \quad  + \p_2\varphi^1 (sx+(1-s)y) \dd s\Bigr).
    \end{align*}

\noindent
Now, the theorems of continuity and limit under the integral sign imply that, for $k\in\{1,2\}$, we have \(\int_0^1 \p_i \varphi^k  (sx+(1-s)y) \dd s\in \C^0_0(\R^2 \times \R^2)\). Thus we obtain the result with $X=(X_1,X_2)$ and \(Q_1(X)=-X_1^2/|X|^2,\ Q_2(X)=-X^2_2/|X|^2,\ Q_3(X)=-X_1X_2/|X|^2\).
\end{proof}

\begin{lemma}\label{lem:limit_measure_diagonal}
    Let \((\mu_N)_N\subset\M_f(\R^d)\) and \(\mu\in\M_f(\R^d)\) be non-atomic. We assume that \(\mu_N \to \mu\) and  \(|\mu_N|\to\nu^d+\nu^a\) narrowly, where \(\nu^d\) is non-atomic and 
    \(\nu^a:= \sum_{j=1}^{ +\infty} \gamma_j \delta_{z_j}\)
    is purely atomic. Let \(Q\in \C^0(\R^2\setminus \{0\})\) be a bounded function. Then the narrow limit $m$ of \( Q(x-y) \dd \mu_N(x) \dd\mu_N(y)\)
    has the form
        \begin{equation*}
        \dd m=Q(x-y) \dd \mu(x) \dd \mu(y)+\sum_{j=1}^{ +\infty} c_j \delta (x-z_j) \delta(y-z_j),
    \end{equation*}
    where \( |c_j|\leq \|Q\|_{L^\infty}|\gamma_j|^2\) for all \(j\in \mathbb{N}^*\). Furthermore if \(\mu_N\geq 0\), for all \(N\) then \(\dd m=Q(x-y) \dd \mu(x) \dd \mu(y)\).
\end{lemma}

\begin{proof}
    By the density of functions of the form \(\sum_{k,l} f_k(x)g_l(y)\) with \(f_k,g_l\in \C^0_0(\R^2)\) in \(\C^0_0(\R^2\times \R^2)\) we obtain the narrow convergence
    \begin{align*}
        |\mu_N| \otimes |\mu_N| & \to \beta:= (\nu^d+\nu^a) \otimes (\nu^d+\nu^a) \\
        & =\nu^d\otimes \nu^d +\nu^d\otimes \nu^a +\nu^a \otimes \nu^d +\nu^a \otimes \nu^a.
    \end{align*}
    It is easy to show, using Fubini's theorem (see for instance  \cite[p.566]{Delort_1991}), that in the situation as above the product measure of the diagonal set $\Delta=\{(x,x')\in\R^2: x=x'\}$ satisfies
    \begin{equation}\label{eq:deltazero}
        \beta\mres\Delta= \nu^a \otimes \nu^a =\sum_{j=1}^{ +\infty} |\gamma_j|^2 \delta (x-z_j) \delta (y-z_j),\qquad \beta\mres\Delta(\cdot):=\beta(\cdot\cap\Delta).
    \end{equation}

    \noindent
    We also have \(|Q(x-y) \mu_N \otimes \mu_N| \leq \|Q\|_{L^\infty} |\mu_N| \otimes |\mu_N|\), since \(Q\) is bounded. Thus, applying Lemma \ref{lem:domination_measure}, we conclude that \(m\) is absolutely continuous with respect to \(\beta\) and \(\|Q\|_{L^\infty} \beta \pm m \geq 0.\) Moreover, we find that the restriction $m\mres\Delta$ of \(m\) to the diagonal is absolutely continuous with respect to \(\beta\mres\Delta\), and hence it is of the form \(\sum_{j=1}^{ +\infty} c_j\delta(x-z_j)\delta(y-z_j)\). But since \(|m\mres\Delta|\leq \beta\mres\Delta\), we find that \(|c_j|\leq \|Q\|_{L^\infty} |\gamma_j|^2\). 
    Next, since \( Q(x-y)\) is continuous outside \(\Delta\), for all \(\varphi \in \C^0_0(\R^2\times \R^2)\), which vanish on \(\Delta\), we have
    \begin{equation*}
        \iint_{\R^2\times \R^2}\varphi(x,y) Q(x-y) \dd \mu_N (x)\dd \mu_N(y) \rightarrow \iint_{\R^2\times \R^2} \varphi(x,y) Q(x-y) \dd \mu(x) \dd \mu(y).
    \end{equation*}
    Let \(\eta \in \C^0_0(\R^2)\) be a function, equal to 1 in a neighborhood of the origin. Then, for all \(\varphi \in \C^0_0(\R^2\times \R^2)\), we have
    \begin{align*}
        \iint_{\R^2 \times \R^2} \varphi(x,y) & Q(x-y) \dd \mu_N(x) \dd \mu_N(y) \\
       & = \iint_{\R^2\times \R^2} \varphi(x,y) Q(x-y) \eta \left( \frac{|x-y|}{\delta}\right)\dd \mu_N(x) \dd \mu_N(y) \\
       & \qquad + \iint_{\R^2\times \R^2} \varphi (x,y) Q(x-y) \left( 1-\eta\left(\frac{|x-y|}{\delta}\right) \right)\dd \mu_N(x) \dd\mu_N(y) \\
        & \xrightarrow[]{n \to  +\infty} \iint_{\R^2\times \R^2} \varphi(x,y) \eta \left( \frac{|x-y|}{\delta}\right) \dd m(x,y)  \\
        & \qquad + \iint_{\R^2 \times \R^2} \varphi(x,y) Q(x-y) \left( 1-\eta\left(\frac{|x-y|}{\delta}\right) \right)\dd \mu(x) \dd \mu(y) \\
        & \xrightarrow[]{\delta \to 0} \iint_{\Delta} \varphi(x,y)\dd m(x,y) 
         + \iint_{\R^2 \times \R^2} \varphi(x,y) Q(x-y) \dd \mu(x) \dd \mu(y),
    \end{align*}
    where we used the dominated convergence theorem to justify the last line.  This concludes the first part of the lemma. Now if we assume \(\mu_N \geq 0\) then \(\mu_N=|\mu_N|\) and \(\nu^a=0\). This implies that all the \(c_j\)'s are zero.
\end{proof}

\noindent
With the above lemmas at hand, we are ready to present the proof of Theorem \ref{th:main_2}.

\begin{proof}[Proof of Theorem \ref{th:main_2}]
    First observe that, for \(\varphi\in\C^\infty_c(\R^2)\), if we multiply equations \eqref{eq:critical_points} by \(\frac{d_i}{M_N}\varphi(x_i)\) and if we sum these equations over \(i=1,\dots,N\)  we obtain
    \begin{multline*}
        \frac{1}{M_N^2}\sum_{i=1}^N \sum_{j\neq i} d_id_j [\nabla g(x_i-x_j)+\nabla F(x_i-x_j)]\cdot \varphi(x_i) +\frac{1}{M_N}\sum_{i=1}^N d_i \nabla V(x_i)\cdot \varphi(x_i)=0. 
    \end{multline*}
    this can be rewritten as
    \begin{multline*}
        \iint_{(\R^2\times \R^2)\setminus \Delta} (\nabla g (x-y)+\nabla F(x-y))\cdot \varphi(x) \dd \mu_N(x) \dd \mu_N(y) \\
        +\int_{\R^2}\nabla V(x)\cdot \varphi(x) \dd \mu_N(x)=0.
    \end{multline*}
    By the anti-symmetry of \(\nabla g\) and \(\nabla F\) the above is equivalent to 
    \begin{multline}\label{eq:passage_limite}
        \frac12 \iint_{(\R^2\times \R^2)\setminus \Delta} (\nabla g (x-y)+\nabla F(x-y))\cdot (\varphi(x)-\varphi(y)) \dd \mu_N(x) \dd \mu_N(y) \\
        +\int_{\R^2}\nabla V(x)\cdot \varphi(x) \dd \mu_N(x)=0.
    \end{multline}
    Note that, under the \(d=2\) assumption, the function \((x,y)\mapsto\nabla g(x-y)\cdot (\varphi(x)-\varphi(y))\) is bounded in \(\R^2\times \R^2\).
    Moreover, we can arbitrarily decide that it is equal to zero on the diagonal \(\Delta\), since $\mu$ is non-atomic, and thus $\Delta$ is  \(\mu \otimes \mu\)-null.    
     With this convention in mind, we aim to pass to the limit in \eqref{eq:passage_limite}.
 
      Let \( Q_i\), \(i=1,2,3\) be as in Lemma \ref{lem:factorization}. Then \( Q_i(x-y) \dd \mu_N(x) \dd \mu_N(y)\) converge narrowly to some \(m_i\in \mathcal{M}(\R^2\times \R^2)\). Therefore, by Lemma \ref{lem:limit_measure_diagonal},  using the notation therein, there exist \(|c_j^{(i)}|\leq |\gamma_j|^2\) for \(i=1,2,3\), \(j\in \mathbb{N}^*\) such that
    \begin{equation*}
        \dd m_i =Q_i(x-y) \dd \mu(x)\dd \mu(y) + \sum_{j=1}^{ +\infty} c_j^{(i)}\delta(x-z_j) \delta (y-z_j).
    \end{equation*}
    Then we pass to the limit in  \eqref{eq:passage_limite} obtaining  
    \begin{multline}\label{eq:tester}
        \iint_{\R^2\times \R^2} \left(\nabla g (x-y) +\nabla F (x-y)\right) \cdot\left( \varphi(x)-\varphi(y)\right) \dd \mu (x) \dd \mu(y) 
        +\int_{\R^2} \varphi(x) \nabla V(x) \dd \mu(x)  \\
        +\sum_{j=1}^{ +\infty} c_j^{(1)} \p_1\varphi^1(z_j)+c_j^{(2)}\p_2\varphi^2(z_j)+c_j^{(3)} (\p_1 \varphi^2+\p_2 \varphi^1)(z_j)=0,
    \end{multline}
    where all \( c^{(i)}_j\in \R\) and \(z_j \in \R^2\) originate from the atomic diagonal part in Lemma \ref{lem:limit_measure_diagonal}.
    A good choice of test functions $\varphi$ reveals that all the coefficients \( c_j^{(i)}\) are equal to zero. Take for instance
    \begin{equation*}
        \varphi^{1,1}_{\delta,j}(x) = (x_1-(z_j)_1) \eta \left(\frac{|x-z_j|^2}{\delta^2} \right), \quad \varphi^{2,1}_{\delta,j}(x)=0,
    \end{equation*}
    where \(\eta\in \C^\infty_c(\R^2)\) is equal to \(1\) in a neighborhood of the origin and vanishes outside \(B_2(0)\). Then
    \[\p_1 \varphi^{1,1}_{\delta,j} = \eta \left(\frac{|x-z_j|^2}{\delta^2} \right)+  \frac{2(x_1-(z_j)_1)^2}{\delta^2}\eta'\left(\frac{|x-z_j|^2}{\delta^2} \right) \xrightarrow[\text{a.e.}]{\delta \to 0} 0, \quad \p_1 \varphi^{1,1}_{\delta,j}(z_j)=1.\]
    Furthermore, \(|\p_1\varphi^{1,1}_{\delta,j}| \leq C \) for some constant \(C>0\) independent of \(\delta\).  Similarly, we have 
    \begin{equation*}
        \p_2 \varphi^{1,1}_{\delta,j}=2 (x_1-(z_j)_1) \frac{x_2-(z_j)_2}{\delta^2}\eta'\left(\frac{|x-z_j|^2}{\delta^2} \right) \xrightarrow[\text{a.e.}]{\delta \to 0} 0,
    \end{equation*}
    and \(|\p_2\varphi^{1,1}_{\delta,j}| \leq C \) for some constant \(C>0\) independent of \(\delta>0.\) By using this information in \eqref{eq:tester}, the dominated convergence theorem shows that \(c_j^{(1)}=0\) for all \(j\in \mathbb{N}\). Similar choices of test functions (e.g. \( \varphi^{1,2}_{\delta,j}=0,\ \varphi^{2,2}_{\delta,j}=(x_2-(z_j)_2)\eta\left(\frac{|x-z_j|^2}{\delta^2} \right)\) and \(\varphi_{\delta_j}^{1,3}=(x_2-(z_j)_2)\eta \left(\frac{|x-z_j|^2}{\delta^2}\right)\)) show that \(c_j^2=c_j^3=0\) for all \(j \in \mathbb{N}\), and hence we conclude the proof.
\end{proof}

In Theorem \ref{th:main_2} we have obtained another critical condition on the limiting measure \(\mu\). In fact, this condition coincides with the condition in Theorem \ref{th:main1.5} thanks to the following proposition which was proved in \cite[p.569]{Delort_1991}
for \(\mu \geq 0\) 
(see also \cite[Lemma 2.1]{Schochet_1995}).

\begin{proposition}\label{prop:equivalence}
Let \(\mu \in \mathcal{M}_f(\R^2)\cap H^{-1}_{{\rm loc}}(\R^2)\), then for all \(\varphi \in \C^\infty_c(\R^2;\R^2)\)
\begin{equation}
    \iint_{\R^2\times \R^2} \nabla g(x-y)\cdot  \left(\varphi (x)-\varphi(y) \right) \dd \mu(x) \dd \mu(y) =\frac{1}{c_d}\int_{\R^2} D \varphi(x) : \lh \nabla g \ast \mu, \nabla g \ast \mu\rh \dd x,
\end{equation}
where \(c_d\) is a dimensional constant  such that \(-\Delta g=c_d \delta_0.\)
\end{proposition}

\begin{proof}
    If we assume that \(\mu \in \C^\infty(\R^2)\) then, by Fubini's theorem we see that 
    \begin{equation}
        \iint_{\R^2\times \R^2}   \nabla g (x-y) \cdot \varphi(x) \mu(x) \mu(y) \dd x \dd y = \int_{\R^2}  (\nabla g \ast \mu)(x) \cdot\varphi(x) \mu(x) \dd x.
    \end{equation}
    hence by symmetry of \(g\), and by using Lemma \ref{lem:chain_rule_stress_tensor}
    \begin{align*}
        \iint_{\R^2 \times \R^2} \nabla g (x-y)\cdot \left(\varphi(x)-\varphi(y) \right) \mu(x) \mu(y) \dd x \dd y &= 2 \int_{\R^2}  (\nabla g \ast \mu)(x) \cdot \varphi(x) \mu(x) \dd x \\
        &=\frac{-1}{c_d} \int_{\R^2} \varphi(x) \cdot \dive \lh \nabla g \ast \mu, \nabla g \ast \mu\rh  \dd x \\
        & =\frac{1}{c_d}\int_{\R^2} D \psi : \lh\nabla g \ast \mu, \nabla g \ast \mu\rh.
    \end{align*}
    Now if \(\mu\) is only in  \(\mathcal{M}_f(\R^2)\cap H^{-1}_{\text{loc}}(\R^2)\) we approximate it by \(\mu_\e:= \rho_\e \ast \mu\), where $\rho_\e$ is the standard mollifier (i.e. \(\rho_\e(x)=\e^{-2}\rho (\e^{-1}x)\) and \(\rho\in \C^\infty_c(\R^2)\), with \(\rho \geq 0\) and \(\int_{\R^2} \rho =1\)). It is well-known, cf.\ \cite[Theorem 2.2]{Ambrosio_Fusco_Pallara_2000}, that \(\mu_\e \in \C^\infty(\R^2)\) and \(\mu_\e \to \mu\) and \(|\mu_\e|\to |\mu|\) weakly-* as finite (signed) Radon measures. Applying Lemma \ref{lem:factorization} and Lemma \ref{lem:limit_measure_diagonal} we find that 
    \begin{multline}\label{eq:equality_1}
        \iint_{\R^2\times \R^2} \nabla g(x-y)\cdot \left( \varphi(x)-\varphi(y)\right)\ \dd \mu_\e (x) \dd \mu_\e(y) \\
        \xrightarrow[]{\e \to +0} \iint_{\R^2\times \R^2}  \nabla g(x-y)\cdot\left( \varphi(x)-\varphi(y)\right) \dd \mu (x) \dd \mu(y).
    \end{multline}
    
    On the other hand, using the Fourier transform, it can be checked that \(\mu_\e \rightarrow \mu\) in \(H^{-1}_{\text{loc}}(\R^2)\). Indeed, let \(\chi \in \C^\infty_c(\R^2)\), then
    \begin{align*}
        \int_{\R^2} \widehat{\chi \mu_\e}(\xi) (1+|\xi|^2)^{-1/2} \dd \xi &= (2\pi)^{-2} \int_{\R^2}  \hat{\mu}_\e \ast \hat{\chi} (\xi)(1+|\xi|^2)^{-1/2} \dd \xi \\
        &=(2\pi)^{-2} \int_{\R^2} \int_{\R^2} \hat{\mu}_\e(\xi-\eta)\hat{\chi}(\eta) (1+|\xi|^2)^{-1/2} \dd \eta \dd \xi \\
        &=(2\pi)^{-2} \int_{\R^2} \int_{\R^2} \hat{\mu}(\xi-\eta) \hat{\rho_\e}(\xi-\eta) \hat{\chi}(\eta) (1+|\xi|^2)^{-1/2} \dd \eta \dd \xi \\
        &=\int_{\R^2} \int_{\R^2} \hat{\mu}(\xi-\eta)\hat{\rho}(\e^{-1}(\xi-\eta))\hat{\chi}(\eta) (1+|\xi|^2)^{-1/2} \dd \eta \dd \xi.
    \end{align*}
    Now, it can be proven that \( \hat{\mu}(\xi-\eta)\hat{\rho}(\e^{-1}(\xi-\eta))\hat{\chi}(\eta) (1+|\xi|^2)^{-1/2}\) converges to \( \hat{\mu}(\xi-\eta)\hat{\chi}(\eta) (1+|\xi|^2)^{-1/2}\) for a.e.\ \( (\xi,\eta)\in \R^2\times \R^2\), since \(\hat{\rho} (\e^{-1} (\xi-\eta))\rightarrow \hat{\rho}(0)=1\) a.e. The function \(\hat{\rho}\) is bounded, and thus, by the dominated convergence theorem,  
    \begin{multline*}
        \int_{\R^2} \int_{\R^2} \hat{\mu}(\xi-\eta)\hat{\rho}(\e^{-1}(\xi-\eta))\hat{\chi}(\eta) (1+|\xi|^2)^{-1/2} \dd \eta \dd \xi \\
        \xrightarrow[]{\e \to 0}\int_{\R^2} \int_{\R^2} \hat{\mu}(\xi-\eta)\hat{\chi}(\eta) (1+|\xi|^2)^{-1/2} \dd \eta \dd \xi,
    \end{multline*}
    and this proves that \(\chi \mu_\e\) converges to \(\chi \mu\) in \(H^{-1}(\R^2)\).
    By elliptic regularity, or better said by an analogue of Proposition \ref{prop:hn}, we have \(\nabla g \ast \mu_\e \rightarrow \nabla g \ast \mu\) in \(L^2_{\text{loc}}(\R^2;\R^2)\). Hence \( \lh \nabla g \ast \mu_\e, \nabla g \ast \mu_\e\rh \rightarrow \lh\nabla g \ast \mu, \nabla g \ast \mu\rh \) in \(L^1_{\text{loc}}(\R^2;\R^{2\times 2})\) and we find that 
    \begin{equation}\label{eq:equality_2}
        \int_{\R^2} D\varphi : \lh\nabla g \ast \mu_\e, \nabla g \ast \mu_\e\rh \dd x
        \xrightarrow{\e \to 0} \int_{\R^2} D\varphi : \lh\nabla g \ast \mu, \nabla g \ast \mu\rh \dd x.
    \end{equation}
    The convergences \eqref{eq:equality_1} and \eqref{eq:equality_2} conclude the proof.
\end{proof}

\section{Mean-field limit of stable critical points}\label{sec:final}
In the final section, we prove Theorem \ref{th:main3}. We assume that \(F,V\) are in \(\C^2(\R^2)\).

\begin{proof}[Proof of Theorem \ref{th:main3}]
Since we assume that each \(\vx_N= (x_1,\dots,x_N)\) satisfies \eqref{eq:critical_points} together with \eqref{eq:Boundedness_Hamiltonian} and Definition \ref{def:stability}, we can multiply \eqref{eq:stability_2} by \(\varphi(x_i)-\varphi(x_j)\) for \(i\neq j\) and \(\varphi \in \C^\infty_c(\R^2,\R^2)\) and sum over \(i\neq j\) to obtain
\begin{equation}
    \frac{1}{M_N^2}\sum_{i=1}^N \sum_{j\neq i} d_i d_j (\varphi(x_i)-\varphi(x_j))^T\left(D^2g(x_i-x_j)+D^2F(x_i-x_j)\right)(\varphi(x_i)-\varphi(x_j)) \geq 0
\end{equation}
    or equivalently
    \begin{equation}
        \iint_{(\R^2\times \R^2) \setminus \Delta} (\varphi(x)-\varphi(y))^T\left(D^2g(x-y)+D^2F(x-y)\right)(\varphi(x)-\varphi(y))\dd \mu_N(x)\dd\mu_N(y) \geq 0.
    \end{equation}
    As for the first order criticality condition, we observe that the kernel \[(\varphi(x)-\varphi(y))^T \left( D^2g(x-y)+D^2F(x-y)\right)(\varphi(x)-\varphi(y))\] belongs to \(\C^0((\R^2\times \R^2)\setminus \Delta)\) and is bounded in all \(\R^2\times \R^2\). Indeed a direct computation shows that 
    \begin{multline*}
       (\varphi(x)-\varphi(y))^T D^2g(x-y)(\varphi(x)-\varphi(y))= -\frac{|\varphi(x)-\varphi(y)|^2}{|x-y|^2} +2\frac{[(\varphi(x)-\varphi(y))\cdot (x-y)]^2}{|x-y|^4}.
    \end{multline*}
    As in the proof of Theorem \ref{th:main_2} we use the fact that $(\mu\otimes\mu)(\Delta) = 0$ to 
    arbitrarily extend the above kernel to the whole $\R^2\times\R^2$. Now, thanks to the assumption that \(\mu_N \geq 0\) for all \(N\), the narrow limit $\mu$ of $\mu_N$ is non-negative too.
   Thus, by Lemma \ref{lem:limit_measure_diagonal} we find that 
    \begin{multline*}
         \iint_{(\R^2\times \R^2) \setminus \Delta} (\varphi(x)-\varphi(y))^T\left( D^2g(x-y)+D^2F(x-y)\right)(\varphi(x)-\varphi(y)) \dd \mu_N(x)\dd\mu_N(y) \\
         \xrightarrow[]{N\to  +\infty} \iint_{(\R^2\times \R^2) \setminus \Delta} (\varphi(x)-\varphi(y))^T \left(D^2g(x-y)+D^2F(x-y)\right) (\varphi(x)-\varphi(y))\dd \mu(x)\dd\mu(y)
    \end{multline*}
    and hence the term in the right-hand side, being the limit of non-negative terms, is also non-negative. This concludes the proof.
\end{proof}

We observe that, in the case of \(V=0\), we recover the expected second-order optimality condition for \(\H\). Note that \(V=0\) satisfies assumption \eqref{A2} in dimension 2. To derive the expected stability condition we have used only \eqref{eq:stability_2} in Definition \eqref{def:stability}. The trouble of using \eqref{eq:stability_1} is that we have to multiply it by \(\varphi(x_i)\), for \(\varphi \in \C^\infty_c(\R^2,\R^2)\), and sum over \(i=1,\dots,N\) to obtain 
\begin{multline*}
    \frac{1}{M_N^2}\sum_{i=1}^N\sum_{j\neq i} d_id_j  \varphi(x_i)^T\left(D^2g(x_i-x_j)+D^2F(x_i-x_j)\right) \varphi(x_i)\\
    +\frac{1}{M_N}\sum_{i=1}^N d_i \varphi(x_i)^T D^2V(x_i) \varphi(x_i)\geq 0.
\end{multline*}
Contrarily to the first order criticality condition we cannot use directly the symmetry of \(D^2g(x-y)+D^2F(x-y)\) to make terms in \(\varphi(x)-\varphi(y)\) appear. This is because the expression is quadratic and mixed terms such as
\begin{equation}
    \iint_{(\R^2\times \R^2)\setminus \Delta}\varphi(x)^T \left(D^2g(x-y)+D^2F(x-y)\right) \varphi(y)\dd \mu_N(x) \dd \mu_N(y)
\end{equation}
appear when trying to symmetrize the kernel.

\section{Conclusion}

In 2D we have shown how to pass to the limit in the first-order and second-order optimality conditions for discrete interaction energies such as $\H_N$ to recover the first- and second-order optimality conditions of continuous interaction energies such as $\H$. This justifies the link between the discrete and the discontinuous energy. We have used two different approaches to obtain the first order criticality condition in the limit and we now discuss the advantages and disadvantages of each approach. The first approach, consists in passing to the limit  in the equation satisfied by the stress-energy tensor of the electric field associated to the empirical measure of critical points. This is a PDE approach and, if we make the parallel with the 2D Euler equations corresponds to a velocity approach. Its main advantage is that it gives a criticality condition without assuming any a priori property  on the limit \(\mu\) of the empirical measures. However it has the disadvantages to rely on the exact form of \(g\) and its link with the Laplacian and the Dirichlet energy. Hence this approach does not seem to apply immediately to interaction \(G\) such that \(|\nabla G(x)|\simeq \frac{1}{|x|}\) but with \(G\) anisotropic. We observe that anisotropic interaction energies are of interest for the study of dislocations and we refer to \cite{Mora_Rondi_Scardia_2019,Carrillo_Mateu_Mora_Rondi_Scardia_Verdera_2021,Mateu_Mora_Rondi_Scardia_Verdera_2023} and references therein. We could not either follow this approach to treat the case of Riesz interaction of the form \(1/|x|^s\) for \(0<s<2\). Indeed in that case the natural approach would be to use an extension in \(\R^3\) and to use the local (but non-constant coefficient) operator associated to the Riesz potential in dimension 3, cf.\ \cite{Serfaty_2020}. However, because our argument does not work in dimension 3, cf.\ lemma 4.4 and the comments below, we were not able to do it. As a side remark we observe that we can treat the case of log interaction in 1D since, with an extension argument it boils down to our study in 2D. The second approach, that we refer to as the vorticity approach by analogy with the 2D Euler equations, has the advantage of working for more general singular interactions since it requires only the kernel \(\nabla g(x-y)\cdot (\varphi(x)-\varphi(y))\) to be bounded. Hence it can be employed to tread anisotropic energies such as the ones appearing in \cite{Mora_Rondi_Scardia_2019}. Another advantage is that we are able to obtain second order criticality condition with this method, assuming the non-negativity of the sequencce of empirical measures. However, as a drawback we need to assume a priori that the limiting measure \(\mu\) is diffuse. It is also not immediate how to deal with the Riesz interactions with this approach.

\bibliographystyle{abbrv}
\bibliography{bib}
\end{document}